

\input amstex
\documentstyle{amsppt}

\input label.def
\input degt.def


\input epsf
\def\picture#1{\epsffile{#1-bb.eps}}
\def\cpic#1{$\vcenter{\hbox{\picture{#1}}}$}

\def\ie{\emph{i.e.}}
\def\eg{\emph{e.g.}}
\def\cf{\emph{cf}}
\def\via{\emph{via}}
\def\etc{\emph{etc}}
\def\viceversa{\emph{vice versa}}

{\catcode`\@11
\gdef\proclaimfont@{\sl}}

\newtheorem{Remark}\Remark\thm\endAmSdef
\conjecture\thm\endproclaim
\def\paragraph{\subsection{}}

\def\dash{\item"\hfill--\hfill"}
\def\Dashes{\widestnumber\item{--}\roster}
\def\endDashes{\endroster}

\loadbold
\def\bA{\bold A}
\def\bD{\bold D}
\def\bE{\bold E}

\def\tA#1{\smash{\tilde\bA#1}}
\def\tD#1{\smash{\tilde\bD#1}}
\def\tE#1{\smash{\tilde\bE#1}}

\let\splus\oplus

\def\F{\Bbb F}

\def\CG#1{\Z_{#1}}
\def\BG#1{\Bbb B_{#1}}

\def\DG#1{\Bbb D_{#1}}
\def\SGSet{\Bbb S}
\def\SG#1{\SGSet_{#1}}

\let\Gs\sigma

\def\SL{\operatorname{\text{\sl SL\/}}}

\def\CC{\Cal C}

\def\tS{\tilde S}

\def\B{\bar B}

\def\PP{\bar P}

\def\Cp#1{\Bbb P^{#1}}
\def\term#1-{$\DG{#1}$-}

\let\Ga\alpha
\let\Gb\beta
\let\Gg\gamma
\def\bGa{\bar\Ga}
\let\Gs\sigma
\let\Gr\rho
\def\1{^{-1}}

\def\ls|#1|{\mathopen|#1\mathclose|}
\let\<\langle
\let\>\rangle
\let\onto\twoheadrightarrow
\def\discr{\operatorname{discr}}
\def\Aut{\operatorname{Aut}}

\def\ord{\operatorname{ord}}
\def\Sk#1{\operatorname{Sk}#1}

\def\bGamma{\bar\Gamma}

\let\bb\delta
\let\bc\beta
\let\MB=M

\def\spec{_{\fam0 sp}}
\def\nonspec{_{\fam0 ns}}

\def\tabstrut{\vrule height9.5pt depth2.5pt}
\def\exstrut{\omit\vrule height2pt\hss\vrule}

\def\inserthyphen{\ifcat\next a-\fi\ignorespaces}
\let\BLACK\bullet
\let\WHITE\circ
\def\CROSS{\vcenter{\hbox{$\scriptstyle\mathord\times$}}}
\let\STAR*
\def\pblack-{$\BLACK$\futurelet\next\inserthyphen}
\def\pwhite-{$\WHITE$\futurelet\next\inserthyphen}
\def\pcross-{$\CROSS$\futurelet\next\inserthyphen}
\def\pstar-{$\STAR$\futurelet\next\inserthyphen}
\def\black{\protect\pblack}
\def\white{\protect\pwhite}
\def\cross{\protect\pcross}
\def\star{\protect\pstar}
\def\NO#1{\mathord\#_{#1}}
\def\nblack{\NO\BLACK}
\def\nwhite{\NO\WHITE}
\def\ncross{\NO\CROSS}
\def\nstar{\NO\STAR}

\def\GAP{{\tt GAP}}
\def\fragment(#1)#2{\ref{fig.e7}(#1)--$#2$}

\def\beginGAP{\bgroup
 \catcode`\^=12\catcode`\#=12\catcode`\_=12
 \obeylines\obeyspaces\eightpoint
 \tt}
\let\endGAP\egroup

\def\Tab{\vtop\bgroup\openup1pt\halign\bgroup$##$\hss\cr}
\def\endTab{\crcr\egroup\egroup}

\topmatter

\author
Alex Degtyarev
\endauthor

\title
Plane sextics via dessins d'enfants
\endtitle

\address
Department of Mathematics,
Bilkent University,
06800 Ankara, Turkey
\endaddress

\email
degt\@fen.bilkent.edu.tr
\endemail

\abstract
We develop a geometric approach to the study of plane sextics with
a triple singular point. As an application, we give an
explicit geometric
description of all irreducible maximal sextics with a
type~$\bold E_7$ singular point and compute their fundamental
groups. All groups found are finite; one of them is
nonabelian.
\endabstract



\endtopmatter

\document

\section{Introduction}

\subsection{Motivation}
The subject of this paper is singular complex plane
projective algebraic curves of degree six
(sextics), considered up to equisingular deformation.
Throughout the paper we assume that all curves involved have
at worst simple singularities.
Formally, the classification of plane sextics can be reduced to a
purely arithmetical problem, see~\cite{JAG},
which can be solved
in many interesting cases,
see, \eg, I.~Shimada's list~\cite{Shimada} of maximal sextics,
or the classification of classical Zariski pairs in
A.~\"Ozg\"uner~\cite{Aysegul}, or the
list of
special
sextics in
A.~Degtyarev~\cite{degt.Oka}; the general
impression is that one can answer any reasonable particular question,
although the complete classification would require
an enormous amount
of work. Furthermore, this arithmetical approach, based on the
theory of $K3$-surfaces, does solve a number of problems
concerning
the geometry of plane sextics, see, \eg, the solution to M.~Oka's
conjecture~\cite{EyralOka.abelian} in Degtyarev~\cite{degt.Oka}
and Tokunaga~\cite{Tokunaga.new}, or
the classification of $Z$-splitting curves in
Shimada~\cite{Shimada.Z}, or the classification
of stable symmetries of
irreducible sextics in Degtyarev~\cite{symmetric}.

However, more subtle questions, such as the computation of the
fundamental group of the complement of a sextic, still remain
unanswered, as they require a much more thorough understanding of
the topology of the curve. A great deal of efforts has been made
lately, see~\cite{degt.8a2}, \cite{degt.e6}, \cite{degt.2a8},
\cite{degt-Oka}, \cite{EyralOka.1}, \cite{EyralOka.2},
\cite{EyralOka.new} (see~\cite{EyralOka.new} for more references)
in order to compute the fundamental groups of relatively few
curves. Each time, the main achievement is discovering a way to
visualize a particular curve
and its braid monodromy; once this is done, computing the
group is a technicality.

Apart from a few curves given by explicit equations, most
approaches to the visualization of plane sextics found in the
literature rely, in one way or another, to an elliptic pencil in
the covering $K3$-surface.
One
such approach
was suggested
in~\cite{symmetric}: one uses a stable symmetry and represents the
sextic as a double covering of an appropriate trigonal curve. In
the present paper, we suggest another approach, which is also
based on the study of trigonal curves; in the long run, we
anticipate to be able to use this correspondence
to handle all sextics with a triple singular
point. Here, we deal with the type~$\bE_7$ singular points and
prepare the background for
types~$\bE_6$ and~$\bE_8$, which are to be the subject of a forthcoming
paper.

It is worth mentioning that, by now, the approach suggested
in~\cite{symmetric} is almost exhausted, at least if one tries to
confine oneself to irreducible maximal trigonal curves: the only
case that has not been considered yet is that of sextics with two
type~$\bE_8$ singular points. In my next paper, it will be shown
that any such sextic has abelian fundamental group.
Jumping a few steps ahead, I can
announce that
{\proclaimfont the only irreducible maximal sextic with a
type~$\bE_8$ singular point and
nonabelian fundamental group has
the set of singularities
$\bE_8\splus\bA_4\splus\bA_3\splus2\bA_2$; its
group
is a semidirect product of its abelianization~$\CG6$ and its
commutant $\SL(2,\F_5)$.}

\subsection{Principal results}
Recall that a plane sextic~$B$ is called \emph{maximal}
(sometimes, \emph{maximizing}), if the total Milnor number
$\mu(B)$ of the singular points of~$B$ takes the maximal possible
value, which is~$19$ (see U.~Persson~\cite{Persson}, where the term was
introduced).
Maximal sextics are projectively rigid; they
are always defined over algebraic number fields.

\theorem\label{th.e7}
Up to projective transformation \rom(equivalently, up to
equisingular deformation\rom),
there are $19$ maximal irreducible plane sextics $B\subset\Cp2$
with simple
singularities only and with at least one
type~$\bE_7$ singular point\rom; they realize $11$ sets of
singularities \rom(see Table~\ref{tab.e7} on
Page~\rom{\pageref{tab.e7}}\rom).
\endtheorem

This theorem is proved in~\ref{s.proof}, where all
sextics are constructed explicitly using trigonal curves.
Alternatively, the statement of the theorem follows from combining
the results of J.-G.~Yang~\cite{Yang} (the existence) and
I.~Shimada~\cite{Shimada} (the enumeration of the sets of
singularities realized by more than one deformation
family). In this respect, it is worth emphasizing that our proof
is purely geometric; although not writing down explicit equations,
we provide a means to completely recover the topology of the pair
$(\Cp2,B)$ and even the topology of the projection
$(\Cp2,B)\to\Cp1$ from the type~$\bE_7$ singular point. (For the
description of the topology of a trigonal curve in terms of its
skeleton, see~\cite{DIK.elliptic} or~\cite{degt.kplets}.)
As an application of this geometric construction, we compute the
fundamental groups $\pi_1(\Cp2\sminus B)$ of
all curves involved and
study their perturbations.

\theorem\label{th.e7.group}
With one exception, the fundamental group $\pi_1(\Cp2\sminus B)$
of a plane sextic $B\subset\Cp2$ as in Theorem~\ref{th.e7} is
abelian. The exception is the \rom(only\rom) sextic with the set
of singularities $\bE_7\splus2\bA_4\splus2\bA_2$\rom; its group is
given by
$$
\multline
G=\bigl\<\Ga_1,\Ga_2,\Ga_3\bigm|
 [\Ga_2,\Ga_3]=[\Ga_i,\Gr^3]=[\Ga_i,\Ga_2^2\Ga_3]=1,
 \ i=1,2,3,\\
\Gr^2\Ga_1=\Ga_2,\
 (\Ga_1\Ga_2)^2\Ga_1=\Ga_2(\Ga_1\Ga_2)^2,\
 (\Ga_1\Ga_3)^2\Ga_1=\Ga_3(\Ga_1\Ga_3)^2\bigr\>,
\endmultline
$$
where $\Gr=\Ga_1\Ga_2\Ga_3$. One can represent
this group~$G$ as a
semi-direct product of its abelianization $\CG6$ and its commutant
$[G,G]\cong\SL(2,\F_{19})$, which is the only perfect group of
order $6840$.
\endtheorem

\theorem\label{th.e7.pert}
For any proper
perturbation~$B'$ of any plane sextic~$B$ as in
Theorem~\ref{th.e7},
the fundamental group $\pi_1(\Cp2\sminus B')$ is
abelian.
\endtheorem

Theorems~\ref{th.e7.group} and~\ref{th.e7.pert} are
proved in~\ref{proof.e7.group} and~\ref{proof.e7.pert},
respectively.

Although we do not treat systematically reducible sextics (the
principal reason being the fact that \GAP~\cite{GAP} does not work
well with infinite groups), the following by-product of our
calculation seems worth mentioning (see~\ref{proof.split} for the
proof).

\proposition\label{th.split}
Let~$B$ be a plane sextic splitting into two irreducible
cubics and having
one of the following five sets of
singularities~$\Sigma$\rom:
$$
\gathered
2\bE_7\splus\bA_5,\qquad
\bE_7\splus\bD_{12},\qquad
\bE_7\splus\bD_5\splus\bA_7,\\
\bE_7\splus\bA_{11}\splus\bA_1,\qquad
\bE_7\splus\bA_9\splus\bA_2\splus\bA_1.
\endgathered
$$
Let $G=\pi_1(\Cp2\sminus B)$. Then, for
$\Sigma=2\bE_7\splus\bA_5$ or $\bE_7\splus\bA_{11}\splus\bA_1$, the
commutant $[G,G]$ is a central subgroup of order~$3$\rom; in the
other three cases, $G$ is abelian.
\endproposition

\corollary\label{th.split.pert}
For any irreducible sextic~$B'$ obtained by a perturbation from a
sextic~$B$ as in Proposition~\ref{th.split}, the group
$\pi_1(\Cp2\sminus B')\cong\CG6$ is abelian.
\endcorollary

Theorems~\ref{th.e7.group} and~\ref{th.e7.pert} and
Corollary~\ref{th.split.pert} provide further evidence to
substantiate my conjecture that the fundamental group of an
irreducible sextic that is not of torus type
(\ie, not given by a polynomial of the form
$p^3+q^2$)
is finite.

Altogether, Theorem~\ref{th.e7.pert}
gives rise to about $250$ new sets of
singularities that are
realized by sextics with abelian fundamental groups,
and Corollary~\ref{th.split.pert} adds about $70$ more.
(Recall that, according to~\cite{degt.8a2},
any induced subgraph of the
combined Dynkin graph of a sextic~$B$ with simple singularities
can be realized by a perturbation of~$B$; in other words, the
singular points of~$B$ can be perturbed independently.) Of special
interest are the eleven (mentioning only the new ones) sets
of singularities listed in Table~\ref{tab.Zariski}. The
corresponding curves are included into the so called
\emph{classical Zariski pairs}, \ie, pairs of plane sextics that
share the same set of singularities but differ by the Alexander
polynomial.
Together with the previously known results, see
Degtyarev~\cite{degt.2a8}
and Eyral, Oka~\cite{EyralOka.abelian}, this makes $46$ out of the
$51$ (see \"Ozg\"uner~\cite{Aysegul}) classical Zariski pairs.
In each
of these $46$ pairs,
the groups of the two curves are $\CG2*\CG3$
and~$\CG6$. (For the curves with nonabelian groups,
see~\cite{degt.2a8} and references there. Note that, formally, all
classical Zariski pairs are known, one of them being in fact a
triple, see~\cite{Aysegul},
but not all curves have been
constructed explicitly, hence
not all fundamental groups
have been computed yet.)

\midinsert
\table\label{tab.Zariski}
`New' classical Zariski pairs
\endtable
\hbox to\hsize{\hss
\Tab
3\bE_6\cr
2\bE_6\splus\bA_5\splus\bA_1\cr
2\bE_6\splus\bA_5\cr
\endTab\hss\Tab
\bE_6\splus\bA_{11}\splus\bA_1\cr
\bE_6\splus\bA_8\splus\bA_2\splus\bA_1\cr
\bE_6\splus2\bA_5\splus\bA_1\cr
\bE_6\splus2\bA_5\cr
\endTab\hss\Tab
\bA_{11}\splus\bA_5\splus\bA_1\cr
\bA_8\splus\bA_5\splus\bA_2\splus2\bA_1\cr
3\bA_5\splus\bA_1\cr
3\bA_5\cr
\endTab
\hss}
\endinsert

Another interesting example is the set of singularities
$3\bA_6$, which is
obtained by a perturbation of $\bE_7\splus2\bA_6$.
The corresponding curve~$B\nonspec$
has a so called \emph{special} counterpart, \ie,
a sextic~$B\spec$ with the set of singularities~$3\bA_6$ and
$\pi_1(\Cp2\sminus B\spec)\cong\CG3\times\DG{14}$, where $\DG{14}$ is
the dihedral group of order~$14$. (The latter group was computed
in Degtyarev, Oka~\cite{degt-Oka}.)
Another, very explicit, construction of a non-special
sextic~$B\nonspec$
with the set of singularities $3\bA_6$
was recently discovered in Eyral, Oka~\cite{EyralOka.new}, where
the group of~$B\nonspec$ was also shown to be abelian. The pair
$(B\spec,B\nonspec)$ constitutes a so called \emph{Alexander
equivalent Zariski pair}: the Alexander polynomials of both curves
are trivial.
We refer
to~\cite{EyralOka.new} for the further discussion of
special sextics and the current state of the subject.

\subsection{Tools and further results}
Our principal tool is to blow up the type~$\bE_7$ point of the
sextic and, after a series of elementary transformations, to
consider the result as a trigonal curve in a ruled rational
surface. We show that maximal sextics
correspond to maximal
trigonal curves, and the latter can be effectively studied using
Grothendieck's \emph{dessins d'enfants} (skeletons in the
terminology of the paper).
In fact, it turns out that a great deal of relevant statements is
already
scattered across~\cite{degt.kplets} and~\cite{symmetric}, and we
merely bring these results together and draw conclusions.

The following intermediate statement seems to be of independent
interest: it gives an estimate on the total Milnor number~$\mu$ of
a non-isotrivial trigonal curve and characterizes maximal curves
as those maximizing~$\mu$. We refer to~\S\ref{S.trigonal} for the
terminology and notation, and to~\ref{proof.mu} for the proof.

\theorem\label{th.mu}
Let $\B$ be a trigonal curve in the Hirzebruch
surface~$\Sigma_k$, and assume that $\B$ is not isotrivial
and that all singularities of~$\B$ are simple.
Then the total Milnor number $\mu(\B)$
of the singular points of~$\B$ is subject to the inequality
$$
\mu(\B)\le5k-2-\#\{\text{\rm unstable fibers of~$\B$}\},
$$
which turns into an equality if and only if $\B$ is maximal.
\endtheorem

It is Theorem~\ref{th.mu} that explains the relation between
maximal sextics and maximal trigonal curves: both maximize the
total Milnor number.

\Remark
The estimate
given by Theorem~\ref{th.mu} does \emph{not}
always hold for
isotrivial trigonal curves,
where
$\mu(\B)$ can be at least as large as $\approx48k/5$.
On the contrary, each non-simple
singular points
reduces the upper bound by~$1$.
\endRemark

Certainly, the crucial property of the sextics in question is the
fact that they have a singular point of multiplicity $d-3$, where
$d$ is the degree of the curve, and our approach is an immediate
generalization of a similar (although much simpler) study of
curves with a singular point of multiplicity $d-2$, see, \eg,
\cite{quintics} and~\cite{groups}. The approach should work
equally well for all sextics with a triple singular point, and
we lay the foundation for a further
development
by completing the
necessary
preliminary calculations for the singular points of type~$\bE_8$
and~$\bE_6$. The precise statements concerning these two types
will appear in a subsequent paper;
the case of a $\bD$-type point will be considered later.


\subsection{Contents of the paper}
In Section~\ref{S.trigonal}, we remind a few basic facts
concerning the
classification of maximal trigonal curves in Hirzebruch surfaces.
The new result proved here is Theorem~\ref{th.mu}.

Section~\ref{S.sextics} introduces the principal tool used in the
paper:
the \emph{trigonal model} of a plane sextic with a
distinguished $\bE$-type singular point. Keeping in mind further
development of the subject, in addition to type~$\bE_7$ appearing
in the principal results (Theorems~\ref{th.e7}--\ref{th.e7.pert}),
in the auxiliary
Sections \ref{S.sextics}--\ref{S.inclusion} we consider as well sextic
with a singular point of type~$\bE_6$ or~$\bE_8$.

In Section~\ref{S.group}, we outline the strategy used to
compute
the fundamental groups, cite a few statements
from~\cite{degt.kplets} concerning the braid monodromy of trigonal
curves, and compute two `universal' relations present in the group
of each curve: the so called relation at infinity and monodromy at
infinity. As a further extension of this analysis of local
canonical forms, in Section~\ref{S.inclusion} we compute the
homomorphism induced by the inclusion of a Milnor ball about an
$\bE$-type singular point.

Theorems~\ref{th.e7} and~\ref{th.e7.group} are proved in
Section~\ref{S.computation}: we enumerate the trigonal models of
sextics with a type~$\bE_7$ singular point by listing their
skeletons (mainly, the problem is reduced to a previously known
classification found in~\cite{symmetric}); then, we use the
skeletons obtained to compute the
fundamental groups.

The concluding Section~\ref{S.perturbations} deals with
Theorem~\ref{th.e7.pert}: first, we (re-)compute the fundamental
groups of the perturbations of a type~$\bE_7$ singular point
(using the same techniques involving trigonal curves and
skeletons), and then we apply these results, the inclusion
homomorphism of Section~\ref{S.inclusion}, and the presentations
obtained in Section~\ref{S.computation} to prove the theorem.

\subsection{Acknowledgements}
This paper was conceived during my
participation in the special semester on Real and Tropical
Algebraic Geometry held at the \emph{Centre Interfacultaire
Bernoulli}, \emph{\'Ecole polytechnique f\'ed\'erale de Lausanne}.
I am thankful to the organizers of the semester and to the
administration of \emph{CIB}.

\section{Trigonal curves}\label{S.trigonal}

In this section, we cite a few results concerning
the classification and properties of maximal trigonal
curves in Hirzebruch surfaces. The principal reference
is~\cite{degt.kplets}.

\subsection{Maximal trigonal curves}\label{s.maximal}
Recall that the \emph{Hirzebruch surface} $\Sigma_k$, $k\ge0$,
is a geometrically ruled
rational surface with an \emph{exceptional section}~$E$
of square~$-k$.
(Sometimes, the fibers of the ruling are referred to as
\emph{vertical lines} in~$\Sigma_k$.)
A \emph{trigonal curve} is a curve $\B\subset\Sigma_k$ disjoint
from~$E$ and
intersecting each generic fiber at three points.

In this paper, {\it we consider trigonal curves with simple
singularities only.}

A \emph{singular fiber} (sometimes also called a
\emph{vertical tangent}) of a trigonal curve $\B\subset\Sigma_k$
is a fiber of the ruling of~$\Sigma_k$ intersecting~$\B$
geometrically at less than three points. Locally,
$\B$ is the
ramification locus of the Weierstra\ss{} model of a Jacobian
elliptic surface, and to describe the
(topological) type of a singular fiber we use
(one of) the standard notation for the singular elliptic fibers,
referring to the extended Dynkin graph of the
corresponding configuration of the exceptional
divisors. The types
are as follows:
\Dashes
\dash
$\tilde\bA_0^*$: a simple vertical tangent;
\dash
$\tilde\bA_0^{**}$: a vertical inflection tangent;
\dash
$\tilde\bA_1^*$: a node of~$\B$ with one of the branches vertical;
\dash
$\tilde\bA_2^*$: a cusp of~$\B$ with vertical tangent;
\dash
$\tilde\bA_p$, $\tilde\bD_q$, $\tilde\bE_6$, $\tilde\bE_7$,
$\tilde\bE_8$: a simple singular point of~$\B$ of the same type
with minimal possible local intersection index with the fiber.
\endDashes
For the relation to Kodaira's classification of singular elliptic
fibers and a few other details, see Table~\ref{tab.j};
further details and references are found in~\cite{degt.kplets}.

\midinsert
\table\label{tab.j}
Types of singular fibers
\endtable
\def\tabrule{\vrule height11pt depth4pt}\def\neg{\!\!\!}
\def\I{\roman{I}}\def\II{\roman{II}}\def\III{\roman{III}}\def\IV{\roman{IV}}

\centerline{\vbox{\offinterlineskip
\halign{\tabrule\quad\hss#\hss\quad\vrule&&\quad\hss$#$\hss\quad\vrule\cr
\noalign{\hrule}
\multispan2\tabrule\hss Type of~$F$\hss\vrule&
 \neg j(F)\neg&\neg\text{Vertex}\neg&\neg\text{Valency}\neg\cr
\noalign{\hrule}
$\tilde\bA_p$ ($\tilde\bD_{p+5}$), $p\ge1$&
 \I_{p+1}\ (\I^*_{p+1})&
 \infty&\text{\cross-}&p+1\cr
\noalign{\hrule}
$\tilde\bA_0^*$ ($\tilde\bD_5$)&\I_1\ (\I^*_1)&
 \infty&\text{\cross-}&1\cr
\noalign{\hrule}
$\tilde\bA_0^{**}$ ($\tilde\bE_6$)&\II\ (\II^*)&
 0&\text{\black-}&1\bmod3\cr
\noalign{\hrule}
$\tilde\bA_1^*$ ($\tilde\bE_7$)&\III\ (\III^*)&
 1&\text{\white-}&1\bmod2\cr
\noalign{\hrule}
$\tilde\bA_2^*$ ($\tilde\bE_8$)&\IV\ (\IV^*)&
 0&\text{\black-}&2\bmod3\cr
\noalign{\hrule}\crcr}}}
\eightpoint\botcaption{\bf Comments}
Fibers of type~$\tilde\bA_0$ (Kodaira's
$\roman{I}_0$) are not singular;
fibers of type~$\tilde\bD_4$
(Kodaira's $\roman{I}_0^*$)
are
not detected by the $j$-invariant.
Fibers of type~$\tilde\bA_0$ or~$\tilde\bD_4$
with complex multiplication of order~$2$
(respectively,~$3$) are over the \white-vertices of valency
$0\bmod2$ (respectively,
over the
\black-vertices of valency $0\bmod3$). The types shown
parenthetically in the table
are obtained from the corresponding $\tA{}$-types
by an elementary transformation, see~\ref{s.elem.tr}; the pairs
are not distinguishable by the $j$-invariant.
\endcaption
\endinsert

The type~$\tA{_0^{**}}$, $\tA{_1^*}$, and~$\tA{_2^*}$
singular fibers of a trigonal curve are called \emph{unstable},
and all other singular fibers are called \emph{stable}.
Informally, a fiber is unstable if its type does not need to
be preserved under equisingular, but not necessarily fiberwise,
deformations of the curve.

A trigonal curve is called \emph{stable} if all its singular
fibers are stable.

The \emph{\rom(functional\rom) $j$-invariant}
$j=j_{\B}\:\Cp1\to\Cp1$ of
a trigonal curve $\B\subset\Sigma_2$ is defined as the analytic
continuation of the function sending a point
$b$ in the base $\Cp1$ of~$\Sigma_2$
to the $j$-invariant (divided
by~$12^3$)
of the elliptic
curve covering the fiber~$F$ over~$b$
and ramified at $F\cap(\B+E)$. The curve~$\B$
is called \emph{isotrivial} if $j_{\B}=\const$. Such curves can
easily be enumerated, see, \eg,~\cite{degt.kplets}.

\definition\label{def.max}
A non-isotrivial trigonal curve~$\B$
is called \emph{maximal} if it has the following properties:
\roster
\item\local{noD4}
$\B$ has no singular fibers of type~$\tD{_4}$;
\item\local{0,1,infty}
$j=j_{\B}$ has no critical values other than~$0$, $1$, and~$\infty$;
\item\local{le3}
each point in the pull-back $j^{-1}(0)$ has ramification index at
most~$3$;
\item\local{le2}
each point in the pull-back $j^{-1}(1)$ has ramification index at
most~$2$.
\endroster
\enddefinition

The maximality of a non-isotrivial trigonal curve $\B\subset\Sigma_2$
can easily be detected by
applying the Riemann--Hurwitz formula to the map
$j_{\B}\:\Cp1\to\Cp1$; it depends only
on the (combinatorial) set of singular fibers of~$\B$,
see~\cite{degt.kplets} for details. An alternative criterion,
based on the total Milnor number,
is given by Theorem~\ref{th.mu} proved in this paper.
The classification of such
curves reduces to a combinatorial problem, see
Theorem~\ref{th.trigonal} below;
a partial classification of
maximal trigonal curves in~$\Sigma_2$ is found
in~\cite{symmetric}. An important property of maximal trigonal
curves is their rigidity, see~\cite{degt.kplets}: any small
deformation of such a curve~$\B$ is isomorphic to~$\B$.

\subsection{Elementary transformations}\label{s.elem.tr}
An \emph{elementary transformation}
(sometimes, \emph{Nagata elementary transformation})
of~$\Sigma_p$
is a birational transformation
$\Sigma_p\dashrightarrow\Sigma_{p+1}$ consisting in blowing up a
point~$P$ in the exceptional section of~$\Sigma_p$ followed by blowing
down (the proper transform of) the fiber~$F$ through~$P$.
(In the sequel, we often omit the reference to
`the proper transform of' when it is understood.)
The inverse transformation
$\Sigma_{p+1}\dashrightarrow\Sigma_{p}$ is also called an
elementary transformation; it consists in blowing up a point~$P'$
\emph{not} in the exceptional section of~$\Sigma_{p+1}$ followed
by blowing down the fiber~$F'$ through~$P'$.

Pick affine charts
$(x_p,y_p)$ in~$\Sigma_p$ and $(x_{p+1},y_{p+1})$ in
$\Sigma_{p+1}$ so that the exceptional sections are given by
$y_p=\infty$ and $y_{p+1}=\infty$, respectively, the fiber~$F$ to
be blown down is given by $x_p=0$, and the image of~$F$ is the
origin $x_{p+1}=y_{p+1}=0$. Then, under the appropriate choice of
$(x_{p+1},y_{p+1})$, the elementary
transformation is the change of coordinates
$$
x_p=x_{p+1},\qquad y_p=y_{p+1}/x_{p+1}.
\eqtag\label{eq.elem.tr}
$$

Let $\B\subset\Sigma_k$ be a \emph{generalized trigonal curve},
\ie, a curve intersecting each generic fiber of the ruling at
three points but possibly \emph{not} disjoint from the exceptional
section~$E$.
Then, by a sequence of elementary transformations, one can resolve
the points of intersection of~$\B$ and~$E$ and obtain a true
trigonal curve $\B'\subset\Sigma_{k'}$, $k'\ge k$, birationally
equivalent to~$\B$. Alternatively, given a trigonal curve
$\B\subset\Sigma_k$ with triple singular points, one can apply a
sequence of elementary transformations to obtain a trigonal curve
$\B'\subset\Sigma_{k'}$, $k'\le k$, birationally
equivalent to~$\B$ and with $\tA{}$~type singular fibers only.

The $j$-invariant~$j_{\B}$, being defined as an analytic
continuation, does not
change under elementary transformations. One can use
this observation to define the
$j$-invariant and all related objects (see next section) for
generalized trigonal curves, as well as for trigonal curves with
triple points, not necessarily simple.

\subsection{Dessins and skeletons}\label{s.dessin}
The concept of
the dessin of a trigonal curve is a modification of Grothendieck's
idea of \emph{dessin d'enfant}; it is due to
S.~Orevkov~\cite{Orevkov}, with a further development
in~\cite{DIK.elliptic} and~\cite{degt.kplets}.

The \emph{dessin}~$\Gamma_{\B}$ of a non-isotrivial trigonal curve
$\B\subset\Sigma_k$ is defined as the planar map
$j_{\B}^{-1}(\Rp1)\subset S^2=\Cp1$,
enhanced with the following decorations: the
pull-backs of~$0$, $1$, and~$\infty$ are called, respectively,
\black--, \white--, and \cross-vertices of~$\Gamma_{\B}$, and the
connected components of the pull-backs of $(0,1)$, $(1,+\infty)$,
and $(-\infty,0)$ are called, respectively, bold, dotted, and
solid edges of~$\Gamma_{\B}$. Clearly, the dessin is invariant
under elementary transformations of the curve and, up to
elementary transformation and isomorphism, the dessin determines
the curve uniquely (see, \eg,~\cite{degt.kplets}; it is essential
that the dessin is considered in the topological sphere~$S^2$; the
analytic structure on~$S^2$
is recovered using the Riemann existence theorem).

The relation between the vertices of the dessin~$\Gamma_{\B}$
and the singular fibers of~$\B$ is shown in Table~\ref{tab.j} (see
also Convention~\ref{valency} concerning the valencies).
The \black-vertices of valency $0\bmod3$ and
\white-vertices of valency $0\bmod2$ are called
\emph{nonsingular}; the other \black-- and \white-vertices are
called \emph{singular}, as they correspond to singular fibers
of the curve.

The \emph{skeleton}~$\Sk_{\B}$ of a trigonal curve~$\B$ is the
planar map obtained from the dessin~$\Gamma_{\B}$ by removing all
\cross-vertices and solid and dotted edges and ignoring all
bivalent \white-vertices. (Thus, $\Sk_{\B}$ is Grothendieck's
\emph{dessin d'enfant} $j_{\B}^{-1}([0,1])$, with the bivalent
pull-backs of~$1$ ignored.) Skeletons are especially useful in the
study of maximal curves. The skeleton~$\Sk_{\B}$ of any
\emph{maximal} curve~$\B$ has the following properties:
\roster
\item\local{Sk.1}
$\Sk_{\B}$ is connected;
\item\local{Sk.2}
each \black-vertex of~$\Sk_{\B}$ has valency~$1$, $2$, or~$3$;
each \white-vertex has valency~$1$ and is connected to a
\black-vertex.
\endroster
Conversely, any planar map $\Sk_{\B}\subset S^2$
satisfying~\loccit{Sk.1},~\loccit{Sk.2} above extends to a unique,
up to orientation preserving diffeomorphism of~$S^2$, dessin of a
\emph{maximal} trigonal curve: one inserts a \white-vertex in the
middle of each edge connecting two \black-vertices, places a
\cross-vertex $u_R$ inside each region~$R$ of~$\Sk_{\B}$, and
connects~$u_R$ by disjoint solid (dotted) edges to all
\black-- (respectively, \white--) vertices in the
boundary~$\partial R$.

According to the following theorem, proved in~\cite{degt.kplets},
skeletons classify maximal trigonal curves. For further
applications of this concept, see~\ref{s.monodromy} below
and~\cite{degt.kplets}.

\theorem\label{th.trigonal}
The fiberwise deformation classes
\rom(equivalently, isomorphism classes\rom) of maximal trigonal curves
$\B\subset\Sigma_k$ with $t$ triple points are in a one-to-one
correspondence with the
orientation preserving diffeomorphism classes of skeletons
\rom(\ie, planar maps $\Sk\subset S^2$ satisfying
conditions~\loccit{Sk.1}, \loccit{Sk.2} above\rom)
satisfying the count given by Corollary~\ref{vertex.count}.
\qed
\endtheorem

\definition[Convention]\label{valency}
It is important to emphasize that the skeleton~$\Sk$ is merely a
convenient way to encode the dessin~$\Gamma$ of a maximal curve;
$\Sk$ is a
subgraph of~$\Gamma$, with some vertices and edges removed and
some vertices ignored. For this reason, in the further exposition
we freely switch between
skeletons and dessins; if only a skeleton is given, we extend it
to the dessin of a maximal curve as explained above.
Convenient in general, this practice may cause a confusion
concerning the valencies of the vertices. To avoid this confusion,
{\em by the valency of a vertex~$v$ we mean one half of the
conventional valency of~$v$ regarded as a vertex of~$\Gamma$\/}. In
other words, we only count the edges of one of the two kinds
present at~$v$. The number thus defined is the conventional
valency of~$v$ in~$\Sk$ (if $v$ \emph{is} a vertex of~$\Sk$);
it also equals the ramification index of~$j_{\B}$ at~$v$.
\enddefinition

\subsection{Markings}\label{s.marking}
Recall that a \emph{marking} at a trivalent \black-vertex~$v$
of a skeleton~$\Sk$ (or dessin~$\Gamma$)
is a counterclockwise ordering $\{e_1,e_2,e_3\}$
of the three edges (respectively,
bold edges)
attached to~$v$.
A marking is uniquely defined by assigning index~$1$ to one of
the three edges.
Given a marking,
the indices of the edges are considered defined
modulo~$3$, so that $e_4=e_1$, $e_5=e_2$, \etc.

A marking at~$v$ defines as well
an ordering $\{e'_1,e'_2,e'_3\}$ of the three
solid edges attached to~$v$: we let~$e'_i$
to be the solid edge opposite to~$e_i$.

A \emph{marking} of the
skeleton~$\Sk$ is a collection of markings at all its trivalent
\black-vertices. Given a marking, one can assign a type $[i,j]$,
$i,j\in\Z_3$, to each edge~$e$ of~$\Sk$ connecting two trivalent
vertices, according to the indices of the two ends of~$e$. A
marking of a skeleton without singular \black-vertices
is called \emph{splitting} if it satisfies the following
two conditions:
\roster
\item
the types of all edges are $[1,1]$, $[2,3]$, or~$[3,2]$;
\item
an edge connecting a \black-vertex~$v$ and a singular
\white-vertex has index~$1$ at~$v$.
\endroster
The following statement is proved in~\cite{degt.kplets}.

\proposition\label{splitting}
A maximal trigonal curve $\B\subset\Sigma_k$
is reducible if and only if its skeleton
has no singular \black-vertices and admits a splitting marking.
\rom(Moreover, each splitting marking defines a component of~$\B$
that is a section of~$\Sigma_k$.\rom)
\qed
\endproposition

\Remark
Proposition~\ref{splitting} is proved by reducing the braid
monodromy, see Subsection~\ref{s.monodromy} below, to the
symmetric group~$\SG3$. A marking at a vertex~$v$
gives rise to a natural ordering
$\{p_1,p_2,p_3\}$ of the three points of the intersection
$F_v\cap\B$, where $F_v$ is the fiber over~$v$, and to a canonical
basis $\{\Ga_1,\Ga_2,\Ga_3\}$ for the fundamental group~$\pi_1$
of the curve,
see~\ref{s.monodromy} or~\cite{degt.kplets} for more details.
For a splitting
marking, the point~$p_1$ over each vertex~$v$ belongs to a
separate component, and in the abelianization
$\pi_1/[\pi_1,\pi_1]$ there is no relation of the form
$\Ga_1=\Ga_2$ or $\Ga_1=\Ga_3$.
The latter observation, combined with the relation at infinity,
see~\ref{s.infty} below, gives one an easy way to find the degree
of the corresponding component.
\endRemark

\subsection{The vertex count}\label{s.count}
Given a
non-isotrivial trigonal curve $\B$,
denote by $\nstar$ the total number of
\star-vertices (where \star-\ stands for
either \black-, or \white-, or
\cross-) in the dessin of~$\B$, let $\nstar(i)$, $i\in\N$, be
the number of \star-vertices of
valency~$i$, and let $\nstar(i\bmod N)$, $i\in\Z_N$, be
the number of \star-vertices of
valency~$i\bmod N$.

Assume that $\B\subset\Sigma_k$ has double singular points only.
Then one has (see~\cite{degt.kplets})
$$
\gather
\deg j_{\B}=\sum_{i>0}i\nblack(i)
 =\sum_{i>0}i\nwhite(i)
 =\sum_{i>0}i\ncross(i),\eqtag\label{eq.deg.j}\\\allowdisplaybreak
6k=\deg j_{\B}
 +2\nblack(1\bmod3)
 +3\nwhite(1\bmod2)
 +4\nblack(2\bmod3),\eqtag\label{eq.6k}\\\allowdisplaybreak
\nblack+\nwhite+\ncross\ge\deg j_{\B}+2,\eqtag\label{eq.R-H}
\endgather
$$
the latter inequality turning into an equality
if and only if $j_{\B}$ has no critical
values other than~$0$, $1$, and~$\infty$, \ie, if $\B$ satisfies
condition~\iref{def.max}{0,1,infty}.

\Remark\label{def.degree}
In~\cite{DIK.elliptic}, the number~$3k$ is called the
\emph{degree} $\deg\Gamma_{\B}$ of the dessin. (The reason is the
fact that generic dessins of degree~$3$ correspond to plane
cubics, regarded as trigonal curves in~$\Sigma_1$.)
We define the degree of a skeleton as the degree of its extension
to the dessin of a maximal curve.
In general,
for a curve $\B\subset\Sigma_k$ with $t$ (simple) triple points,
one has $\deg\Gamma_{\B}=3(k-t)$, \cf\.~\ref{s.elem.tr}. In this
notation, one can replace~$6k$ with
$2\deg\Gamma_{\B}$ in~\eqref{eq.6k}
and lift the assumption that $\B$ should have
double singular points only.
\endRemark

Next statement is an immediate consequence
of~\eqref{eq.deg.j}--\eqref{eq.R-H}.

\corollary\label{vertex.count}
Let $\B\subset\Sigma_k$ be a maximal trigonal curve. Then the
numbers of vertices in its skeleton are subject to the identity
$$
\nblack+\nwhite(1)+\nblack(2)=2(k-t),
$$
where $t$ is the number of triple singular points of~$\B$.
\endcorollary

\proof
By a sequence of elementary transformations one can convert~$\B$
to another maximal trigonal curve $\B'\subset\Sigma_{k-t}$,
see~\ref{s.elem.tr}, which has double singular points only and
has the same dessin as~$B$. Then,
it suffices to substitute to~\eqref{eq.6k} the first expression
for $\deg j_{\B}$ from~\eqref{eq.deg.j}, collect (partially for $i=2$)
the terms with $\nblack(i)$, $i=1,2,3$, and divide by~$3$.
\endproof

\subsection{Proof of Theorem~\ref{th.mu}}\label{proof.mu}
First, assume that $\B$ has double singular points only. Then
$$
\mu(\B)=
 \sum_{i>0}(i-1)\ncross(i)
 +\nwhite(1\bmod2)
 +2\nblack(2\bmod3),
$$
see Table~\ref{tab.j}.
Substituting to~\eqref{eq.R-H} the
third expression for $\deg j_{\B}$ from~\eqref{eq.deg.j}, one
obtains the estimate
$$
\mu(\B)\le\nblack+\nwhite
 +\nwhite(1\bmod2)
 +2\nblack(2\bmod3)-2,
\eqtag\label{eq.ncross}
$$
which is sharp if and only if $\B$ satisfies
condition~\iref{def.max}{0,1,infty}.
Substituting to \eqref{eq.6k}
the first two expressions for $\deg j_{\B}$
and replacing the valencies of \black-- and \white-vertices
with their residues modulo~$3$
(in the range $\{1,2,3\}$) and~$2$ (in the range $\{1,2\}$),
respectively,
one obtains the inequalities
$$
\gather
2k\ge\nblack
 +\nwhite(1\bmod2)
 +\nblack(2\bmod3),\eqtag\label{eq.2k}\\
3k\ge\nwhite
 +\nblack(1\bmod3)
 +\nwhite(1\bmod2)
 +2\nblack(2\bmod3),\eqtag\label{eq.3k}
\endgather
$$
which turn into equalities if and only if the dessin of~$\B$ has
no \black-vertices of valency greater than~$3$
(for~\eqref{eq.2k}) or no
\white-vertices of valency greater than~$2$
(for~\eqref{eq.3k}), \ie, if $\B$
satisfies conditions \iref{def.max}{le3} or~\ditto{le2},
respectively. Combining this with~\eqref{eq.ncross}
and taking into account the
fact that $\nblack(1\bmod3)+\nwhite(1\bmod2)+\nblack(2\bmod3)$ is
the number of unstable fibers of~$\B$, one obtains the desired
inequality. The equality holds if and only if $\B$ is maximal, as
the remaining condition~\iref{def.max}{noD4} holds automatically.

If $\B$ has triple points, one can remove them one by one and use
induction. Let~$\PP$ be a triple point of~$\B$ and let
$\B'\subset\Sigma_{k-1}$ be the curve obtained from~$\B$ by the
inverse elementary transformation centered at~$\PP$. If $\PP$ is
of type~$\bD_4$ (and hence $\B$ is not maximal), then
$\mu(\B)=\mu(\B')+4$ and the inequality for~$\B'$ turns into a
strict inequality for~$\B$. In all other cases, $\B$ and~$\B'$ are
or are not maximal simultaneously. If $\PP$ is of type~$\bD_p$,
$p\ge5$, then $\mu(\B)=\mu(\B')+5$ and the unstable fibers
of~$\B'$ are in a one-to-one correspondence with those of~$\B$.
If $\PP$ is of type~$\bE_6$,
$\bE_7$, or~$\bE_8$, then $\mu(\B)=\mu(\B')+6$ and, on the other
hand, $\B'$ has one extra unstable fiber compared to~$\B$ (of
type~$\tA{_0^{**}}$, $\tA{_1^*}$, or~$\tA{_2^*}$, respectively).
In both cases, the defect in the inequality for~$\B'$ is the same
as the
defect in the resulting inequality for~$\B$, and the statement
follows.
\qed

\section{Plane sextics}\label{S.sextics}

The bulk of this section deals with the trigonal models of plane
sextics, which are the principal tool in both the classification
and the computation of the fundamental group. To facilitate the
classification, we also prove Proposition~\ref{sing.points},
restricting the sets of singularities of irreducible
maximal sextics.

\subsection{Irreducible maximal sextics}
Recall that a plane sextic~$B$ is called \emph{maximal} if its
total Milnor number $\mu(B)$ takes the maximal possible
value~$19$.

\proposition\label{sing.points}
An irreducible maximal plane sextic cannot have a singular
point of type~$\bD_{2k}$, $k\ge2$ or more than one singular point
from the following list\rom:
$\bA_{2k+1}$, $k\ge0$, $\bD_{2k+1}$, $k\ge2$, or~$\bE_7$.
\endproposition

\proof
Formally, one can derive the statement from Yang's
list~\cite{Yang} of the sets of singularities realized by
irreducible maximal sextics. For a more conceptual proof,
consider the double covering of the plane ramified at the sextic
and denote by~$X$ its minimal resolution of singularities. It is a
$K3$-surface. Let $L=H_2(X)$, let $\Sigma\subset L$ be the
sublattice spanned by the classes of the exceptional divisors
contracted by the projection $X\to\Cp2$, and let
$S=\Sigma\oplus\<h\>$, where $h$ is the class realized by the
pull-back of a generic line. One has $h^2=2$ and $\Sigma$ is the
direct sum of (negative definite)
irreducible root systems of the same type ($\bA$--$\bD$--$\bE$)
as the singular points of the sextic. Let $\tS\supset S$ be the
primitive hull of~$S$ in~$L$. As shown in~\cite{JAG}, the quotient
$\tS/S$ is free of $2$-torsion. Hence, the $2$-torsion of the
discriminant groups $\discr S=\discr\Sigma\oplus\Z_2$
and $\discr\tS$ coincide. On the
other hand, $\discr\tS\cong\discr S^\perp$ and, since
$\rank S^\perp=2$, the $2$-torsion of $\discr\Sigma$ must be
a cyclic group.
\endproof

\subsection{Trigonal models: the statements}\label{s.models}
In Propositions~\ref{1-1.e8}, \ref{1-1.e7} and~\ref{1-1.e6},
we introduce
certain trigonal curves birationally equivalent to plane sextics;
these curves will be called the \emph{trigonal model}
of the corresponding sextics. Proofs are given
in~\ref{proof.1-1.e8}--\ref{proof.1-1.e6} below.

\proposition\label{1-1.e8}
There is a natural
bijection~$\phi$, invariant under equisingular deformations,
between the following two sets\rom:
\roster
\item\local{e8.sextic}
plane sextics~$B$ with a distinguished
type~$\bE_8$ singular point~$P$, and
\item\local{e8.trigonal}
trigonal curves $\B\subset\Sigma_3$ with a
distinguished type~$\tA{_1^*}$ singular fiber~$F$.
\endroster
A sextic~$B$ is irreducible if and only if so is $\B=\phi(B)$ and,
with one exception, $B$ is maximal if and only if $\B$ is
maximal and has no unstable fibers other than~$F$.
\rom(The exception is the reducible sextic~$B$ with the set of
singularities $\bE_8\splus\bE_7\splus\bD_4$\rom; in this case,
$\phi(B)$ is isotrivial.\rom)
Furthermore, for each pair $B$, $\B=\phi(B)$, there is a
diffeomorphism
$$
\Cp2\sminus(B\cup L)\cong\Sigma_3\sminus(\B\cup E\cup F),
$$
where $L$ is the line tangent to~$B$ at~$P$
and $E$ is the exceptional section.
\endproposition

\proposition\label{1-1.e7}
There is a natural
bijection~$\phi$, invariant under equisingular deformations,
between Zariski open \rom(in each equisingular
stratum\rom) subsets of the following two sets\rom:
\roster
\item\local{e7.sextic}
plane sextics~$B$ with a distinguished
type~$\bE_7$ singular point~$P$ and without linear components
through~$P$, and
\item\local{e7.trigonal}
trigonal curves $\B\subset\Sigma_3$ with a distinguished
type~$\tA{_1}$ singular fiber $F$ and a distinguished
branch at the corresponding type~$\bA_1$ singular point of~$\B$.
\endroster
A sextic~$B$ is irreducible if and only if so is $\B=\phi(B)$, and
$B$ is maximal if and only if $\B$ is
maximal and stable.
Furthermore, for each pair $B$, $\B=\phi(B)$, there is a
diffeomorphism
$$
\Cp2\sminus(B\cup L)\cong\Sigma_3\sminus(\B\cup E\cup F),
$$
where $L$ is the line tangent to~$B$ at~$P$
and $E$ is the exceptional section.
\endproposition

\Remark
Thus, one should expect that, in many cases,
a maximal stable pair $(\B,F)$
as in Proposition~\iref{1-1.e7}{e7.trigonal} would correspond to
\emph{two} deformation classes of sextics. This is indeed the
case, see the sets of singularities marked with a $^*$ in
Table~\ref{tab.e7} in Subsection~\ref{s.proof}.
Arithmetically, this phenomenon is
probably due to the fact that the discriminant group $\discr S$
has two essentially different copies of $\<\frac12\>$, namely,
those coming
from $\discr\bE_7$ and from $\discr\<h\>$. For details,
see~\cite{JAG}.
\endRemark

\proposition\label{1-1.e6}
There is a natural
bijection~$\phi$, invariant under equisingular deformations,
between Zariski open \rom(in each equisingular
stratum\rom) subsets of the following two sets\rom:
\roster
\item\local{e6.sextic}
plane sextics~$B$ with a distinguished
type~$\bE_6$ singular point~$P$, and
\item
trigonal curves $\B\subset\Sigma_4$ with a distinguished
type~$\tA{_5}$ singular fiber $F$.
\endroster
A sextic~$B$ is irreducible if and only if so is $\B=\phi(B)$, and
$B$ is maximal if and only if $\B$ is
maximal and stable.
Furthermore, for each pair $B$, $\B=\phi(B)$, there is a
diffeomorphism
$$
\Cp2\sminus(B\cup L)\cong\Sigma_4\sminus(\B\cup E\cup F),
$$
where $L$ is the line tangent to~$B$ at~$P$
and $E$ is the exceptional section.
\endproposition

\Remark
There are statements similar to Propositions~\ref{1-1.e8},
\ref{1-1.e7}, and \ref{1-1.e6} for sextics with a distinguished
type~$\bD$ singular point. In this case, one would need
to keep track of
two (three in the case~$\bD_4$) singular fibers of~$\B$.
\endRemark

\Remark
Informally, the relation between maximal sextics and maximal
trigonal curves follows from the fact that both objects are rigid,
\ie, curves are isomorphic to their small equisingular
deformations. Formal proofs are given below.
\endRemark


\subsection{Proof of Proposition~\ref{1-1.e8}}\label{proof.1-1.e8}
The bijection~$\phi$ and the diffeomorphism in the statement are
given by a birational transformation
$\Cp2\dashrightarrow\Sigma_3$, so that~$\B$ is the proper
transform of~$B$: one blows up the distinguished
type~$\bE_8$ point~$P$ to get a
generalized trigonal
curve $B'\subset\Sigma_1$ with a cusp tangent to the exceptional
section of~$\Sigma_1$ (the exceptional divisor of the blow-up),
and then
one applies two elementary transformations
to make the curve disjoint
from the exceptional section, see~\ref{s.elem.tr}.

Pick affine coordinates $(u,v)$ in~$\Cp2$ centered at the
distinguished singular point~$P$ and with the $v$-axis along the
line~$L$ in the statement. By the B\'ezout theorem, $B$
intersects~$L$
at one more point $v=a\ne0$. Hence, up to higher order
terms, $B$ is given by a polynomial of the form
$$
(u^3-v^5)(v-a).
\eqtag\label{eq.curve.e8}
$$
In appropriate affine coordinates $(x,y)=(x_3,y_3)$ in~$\Sigma_3$,
\cf\.~\ref{s.elem.tr}, the transformation is given by the
coordinate change
$$
u=x^3\!/y,\qquad v=x^2\!/y,
\eqtag\label{eq.tr.e8}
$$
(in particular, it restricts to a biholomorphism
$\Cp2\sminus L\to\Sigma_3\sminus(E\cup F)$),
and the proper transform~$\B$ of~$B$ is given by
$$
(y^2-x)(x^2-ay).
\eqtag\label{eq.image.e8}
$$
One can see that $F=\{x=0\}$ is a type~$\tA{_1^*}$
singular fiber of~$\B$.

The construction is obviously invertible: given a trigonal curve
$\B\subset\Sigma_3$ with a type~$\tA{_1^*}$ singular fiber~$F$, one
can apply two elementary transformations centered at (the
transform of) the branch of~$\B$ not tangent to~$F$ and then blow
down the exceptional section of the resulting Hirzebruch
surface~$\Sigma_1$; the result is a sextic with a type~$\bE_8$
singular point.

Since $\B$ and~$B$ are the proper transforms of each other, it is
immediate that $\B$ is irreducible if and only if so is~$B$. The
assertion on maximal curves follows from Theorem~\ref{th.mu}.
Indeed, the total Milnor numbers of~$B$ and~$\B$
are related \via\ $\mu(\B)=\mu(B)-7$: the
singular points of~$B$ are in a one-to-one correspondence with
those of~$\B$, which are of the same type, except that the
type~$\bE_8$ point~$P$ corresponds to the type~$\bA_1$ singular
point of~$\B$ in the fiber~$F$. Hence, $B$ is maximal if and only
if $\mu(\B)=12$. If $\B$ is not isotrivial (the isotrivial case
is treated in the next paragraph), then,
taking into account the fact that
$\B$ does have an unstable fiber~$F$,
Theorem~\ref{th.mu} implies that the
latter equality holds if and only if $\B$ is maximal and has no
other unstable fibers.

If $\B$ is isotrivial, then
$j_{\B}\equiv1$
(as $j_{F}=1$, see Table~\ref{tab.j})
and, in appropriate affine coordinates $(x,y)$, the
curve is given by the Weierstra{\ss} equation of the form
$y^3-yxp(x)=0$, $\deg p=5$ and $p(0)\ne0$. (We assume that $x=0$
is the distinguished type~$\tA{_1^*}$ fiber~$F$.)
Such a curve has singular points of
types $\bA_1$, $\bD_4$, and~$\bE_7$ (corresponding, respectively,
to simple, double, and triple roots of the equation $xp(x)=0$),
and the only such set of singularities with the total Milnor
number~$12$ is $\bE_7\oplus\bD_4\oplus\bA_1$, the
type~$\bA_1$ point being located in~$F$. (Note that this set of
singularities cannot be realized by a maximal non-isotrivial
curve, see~\iref{def.max}{noD4}.)
All other statements are straightforward.
\qed

\subsection{Proof of Proposition~\ref{1-1.e7}}\label{proof.1-1.e7}
As in the previous proof,
the bijection~$\phi$ and the diffeomorphism
are given by a birational transformation
$\Cp2\dashrightarrow\Sigma_3$:
one blows up the distinguished
type~$\bE_7$ point~$P$ to get a generalized trigonal
curve $B'\subset\Sigma_1$ that has a node with one of the branches
tangent to the exceptional
section of~$\Sigma_1$; then,
two elementary transformations centered at this branch
produce a trigonal curve $\B\subset\Sigma_3$, see~\ref{s.elem.tr}.
In appropriate affine
coordinates $(u,v)$ in~$\Cp2$ and $(x,y)$ in~$\Sigma_3$, such that
$L$ is the $v$-axis and $F$ is the $y$-axis,
the
transformation is given by~\eqref{eq.tr.e8}. Up to higher order
terms, the defining
polynomial of a typical (see below) sextic~$B$ as in the statement
has the form
$$
(u^2-v^3)(u-bv^2)(v-a),
\eqtag\label{eq.curve.e7}
$$
$a,b=\const$ (the smooth branch of~$B$
at~$P$ is tangent to~$L$ and $B$ intersects~$L$
at one more point $v=a$), and the transform~$\B$ of~$B$ is given
by the polynomial
$$
(y-1)(y-bx)(x^2-ay).
\eqtag\label{eq.image.e7}
$$
One can see that
$F=\{x=0\}$ is a type~$\tA{_1}$ singular fiber of~$\B$ and
$\PP=(0,0)$ is a type~$\bA_1$ singular point. The branch $x^2=ay$
of~$\B$ at~$\PP$ is the transform of the `separate' branch $v=a$
of~$B$; thus, it is distinguished.

The inverse construction consists in applying two elementary
transformations centered at (the transform of) the distinguished
branch of~$\B$ at~$\PP$, followed by blowing down the exceptional
section of the resulting Hirzebruch surface~$\Sigma_1$.

The assertion on the correspondence between irreducible and
maximal curves is proved similar to~\ref{proof.1-1.e8}. This time,
one has $\mu(\B)=\mu(B)-6$ (the type~$\bE_7$ singular point~$P$ is
replaced with the type~$\bA_1$ singular point~$\PP$); hence, $B$
is maximal if and only if $\mu(\B)=13$, and Theorem~\ref{th.mu}
implies that the latter identity holds if and only if $\B$ is
maximal and stable. Note that $\B$ cannot be
isotrivial, as it has a singular fiber~$F$ of type~$\tA{_1}$
and $j_{\B}(F)=\infty$, see Table~\ref{tab.j}.

It remains to show that $\phi$ is defined on a Zariski open subset
of each equisingular stratum. The only extra degeneration that a
sextic~$B$ within a given stratum may have is that
the smooth branch of~$B$ at~$P$ may become
inflection tangent to~$L$. Then,
the singular fiber~$F$ of~$\B$ is of type~$\tA{_2}$
rather than~$\tA{_1}$. However, from the theory of trigonal curves it
follows that such a fiber can be perturbed to a fiber of
type~$\tA{_1}$ and a close fiber of type~$\tA{_0^*}$; this
perturbation can obviously be followed by a one-parameter family
of inverse birational transformations
$\Sigma_3\dashrightarrow\Cp2$ and hence by an equisingular
deformation of sextics.
\qed

\Remark\label{rem.Shimada}
At the end of the proof, we essentially show, using deformations
of trigonal curves, that a line inflection tangent to the smooth
branch of a type~$\bE_7$ singular point of a plane sextic cannot
be stable under equisingular deformations of the sextic.
Alternatively, one can argue that, if such a line were stable, it
would be a \emph{$Z$-splitting curve} in the sense of
Shimada~\cite{Shimada.Z}, and refer to the classification of
$Z$-splitting curves found in~\cite{Shimada.Z}. A similar
observation applies as well to the end of the proof
in~\ref{proof.1-1.e6}.
\endRemark

\subsection{Proof of Proposition~\ref{1-1.e6}}\label{proof.1-1.e6}
The bijection~$\phi$ and the diffeomorphism are
given by a birational transformation
$\Cp2\dashrightarrow\Sigma_4$: one blows up the distinguished
point~$P$ to obtain a generalized trigonal
curve $B'\subset\Sigma_1$ with a branch
inflection tangent to the exceptional section of~$\Sigma_1$, and
then applies three elementary transformations centered at (the
transforms of) this branch to make the curve disjoint from the
exceptional section, see~\ref{s.elem.tr}.

In appropriate affine coordinates $(u,v)$ in~$\Cp2$ and $(x,y)$
in~$\Sigma_3$, such that $L$ is the $v$-axis and $F$ is the
$y$-axis, the transformation is given by
$$
u=x^4\!/y,\qquad v=x^3\!/y.
\eqtag\label{eq.tr.e6}
$$
A typical (see below) sextic~$B$ as
in the statement intersects~$L$ at two other points
$v=a_1$, $v=a_2$,
$a_1\ne a_2$, $a_1,a_2\ne0$. Hence,
up to higher order terms, its defining
polynomial
has the form
$$
(u^3-v^4)(v-a_1)(v-a_2),
\eqtag\label{eq.curve.e6}
$$
and its transform $\B\subset\Sigma_4$ is given
by the polynomial
$$
(y-1)(x^3-a_1y)(x^3-a_2y).
\eqtag\label{eq.image.e6}
$$
One can see that
$F=\{x=0\}$ is a type~$\tA{_5}$ singular fiber of~$\B$. The inverse
correspondence is given by three elementary transformations
centered at (the transforms of) the type~$\bA_5$ singular point,
followed by blowing down the exceptional section of the
resulting Hirzebruch surface~$\Sigma_1$.

The correspondence between irreducible and
maximal curves is established as above:
one has $\mu(\B)=\mu(B)-1$ (the type~$\bE_6$ singular point~$P$ is
replaced with a type~$\bA_5$ singular point of~$\B$); hence, $B$
is maximal if and only if $\mu(\B)=18$, and Theorem~\ref{th.mu}
implies that the latter identity holds if and only if $\B$ is
maximal and stable. Note that $\B$ cannot be
isotrivial, as
$j_{\B}(F)=\infty$, see Table~\ref{tab.j}.

The only extra degeneration
that a sextic~$B$ may have within a given equisingular stratum is
that it may become tangent to~$L$ or one of its type
$\bA_p$, $p\ge1$, singular points may slide into~$L$.
Then, the
fiber~$F$ of the transform~$\B$
is of type~$\tA{_6}$ or $\tA{_{6+p}}$,
respectively.
Such a fiber can be perturbed to a fiber of
type~$\tA{_5}$ and a close fiber of type~$\tA{_0^*}$
or~$\tA{_p}$, respectively, and this
perturbation is followed by an equisingular
deformation of sextics. Thus, the bijection~$\phi$ is well defined
on a Zariski open subset of each stratum.
(An alternative proof of this fact is
explained in Remark~\ref{rem.Shimada}.)
\qed

\section{The fundamental group}\label{S.group}

We outline the strategy used to compute the fundamental groups,
explain how the braid monodromy can be found, and compute a few
`universal' relations, present in the group of any curve in
question.

\subsection{The strategy}\label{s.strategy}
In this section,
we consider a plane sextic~$B$ with a distinguished type~$\bE$
singular point~$P$ and
use Propositions~\ref{1-1.e8}, \ref{1-1.e7}, and~\ref{1-1.e6}
to transform it to a trigonal curve
$\B\subset\Sigma_k$, $k=3$ or~$4$,
with a distinguished singular fiber~$F$. (We
assume~$B$ generic in its equisingular deformation class.) In each
case, $\B$ has a unique singular point in~$F$; it will be denoted
by~$\PP$.
The
above cited propositions give a diffeomorphism
$$
\Cp2\sminus(B\cup L)\cong\Sigma_k\sminus(\B\cup E\cup F),
$$
where $L$ is the line tangent to~$B$ at~$P$. Hence, there is an
isomorphism
$$
\pi_1(\Cp2\sminus(B\cup L))\cong
 \pi_1(\Sigma_k\sminus(\B\cup E\cup F))
$$
of the fundamental groups. According to
E.~R.~van~Kampen~\cite{vanKampen}
(see also T.~Fujita \cite{Fujita}),
the passage from
$\pi_1(\Cp2\sminus(B\cup L))$ to $\pi_1(\Cp2\sminus B)$ results in
adding an extra relation, which can be represented in the form
$[\partial\Gamma]=1$, where $\Gamma\subset\Cp2$ is a small
holomorphic disk transversal to~$L$ and disjoint from~$B$,
and $[\partial\Gamma]\subset\pi_1(\Cp2\sminus(B\cup L))$ is the
class of the boundary of~$\Gamma$ (more precisely, its conjugacy
class).
Denoting by~$\bGamma$ the image of~$\Gamma$ in~$\Sigma_k$, one has
$$
\pi_1(\Cp2\sminus B)\cong
 \pi_1(\Sigma_k\sminus(\B\cup E\cup F))/[\partial\bGamma].
\eqtag\label{eq.bGamma}
$$
The relation $[\partial\bGamma]=1$ is called the \emph{relation at
infinity}; the bulk of this section deals with computing this
relation.

The group $\pi_1(\Sigma_k\sminus(\B\cup E\cup F))$ is computed
using the classical Zariski--van Kampen method~\cite{vanKampen}.
Pick some coordinates $(x',y')$ in the affine chart
$\Sigma_k\sminus(E\cup F)$. For the further exposition, it is
convenient to take
$$
x'=1/x,\qquad y'=y/x^k,
\eqtag\label{eq.coord}
$$
where $(x,y)$ are the coordinates about~$F$ introduced
in~\ref{proof.1-1.e8}--\ref{proof.1-1.e6}. Let $F_1,\ldots,F_r$ be
all singular fibers of~$\B$ other than~$F$, and let~$F_0$ be a
nonsingular fiber.
Pick a closed topological
disk~$\Delta$ in the $x$-axis containing all~$F_j$,
$j=0,\ldots,r$, in its interior
and let
$\Delta^\circ=\Delta\sminus\{F_1,\ldots,F_r\}$.
(We identify fibers with their projections to the
base of the ruling.) Pick a topological section
$s\:\Delta\to\Sigma_k$ proper in the sense of~\cite{degt.e6}. (For
all practical purposes, it suffices to consider a constant section
$y=c$, where $c$ is a constant, $\ls|c|\gg0$. In~\cite{degt.e6},
one can find a more formal exposition; in particular, it is shown
there that the result does not depend on the choice of a proper
section.
A similar approach is found in~\cite{Artal}.)
Let $\{\Ga_1,\Ga_2,\Ga_3\}$ be a basis for the free group
$\pi_F:=\pi_1(F_0\sminus(\B\cup E),s(F_0))$, and let
$\Gg_1,\ldots,\Gg_r$ be a basis for the group
$\pi_1(\Delta^\circ,F_0)$. Dragging the nonsingular fiber along a
loop~$\Gg_j$, $j=1,\ldots,r$ and keeping the base point in~$s$,
one obtains an automorphism $m_j\in\Aut\pi_F$, which is called the
\emph{braid monodromy} along~$\Gg_j$.
(Since the reference section is proper, this automorphism is
indeed a braid.)
In this notation, the
Zariski--van Kampen theorem states that
$$
\pi_1(\Sigma_k\sminus(\B\cup E\cup F))=
 \bigl<\Ga_1,\Ga_2,\Ga_3\bigm|
 \text{$m_j=\id$, $j=1,\ldots,r$}\bigr>,
\eqtag\label{eq.vanKampen}
$$
where each \emph{braid relation} $m_j=\id$, $j=1,\ldots,r$,
should be understood as
the triple of relations $m_j(\Ga_i)=\Ga_i$, $i=1,2,3$.

\Remark
Since each $m_j$ is a braid and thus preserves $\Ga_1\Ga_2\Ga_3$,
it would suffice to keep the relations $m_j(\Ga_1)=\Ga_1$ and
$m_j(\Ga_2)=\Ga_2$ only. Note however that, in a more advanced
setting, the braid monodromy does not necessarily take values in
the braid group, and all three relations should be kept. Besides,
following the principle `the more relations the better', often it
is more convenient to restate the braid relations in the form
$m_j(\Ga)=\Ga$ for each $\Ga\in\<\Ga_1,\Ga_2,\Ga_3\>$.
\endRemark

We will also consider the \emph{monodromy at infinity}~$m_\infty$,
\ie, the braid monodromy along the loop~$\partial\Delta$
(assuming that the base point~$F_0$ is chosen in $\partial\Delta$).

\proposition\label{presentation}
Let $\{\Gg_1,\ldots,\Gg_r\}$ be a
basis for the
free group $\pi_1(\Delta^\circ,F_0)$ such that
$\Gg_1\ldots\Gg_r=[\partial\Delta]$. Then the group
$\pi_1(\Cp2\sminus B)$ has a presentation of the form
$$
\bigl<\Ga_1,\Ga_2,\Ga_3\bigm|
 \text{$m_j=\id$, $j=1,\ldots,r$, $m_\infty=\id$,
 $[\partial\bGamma]=1$}\bigr>.
$$
Furthermore, in the presence of the last two relations,
\rom(any\rom) one of the first $r$ braid relations $m_j=\id$ can
be omitted.
\endproposition

\proof
The presentation is given by~\eqref{eq.bGamma}
and~\eqref{eq.vanKampen}; the relation $m_\infty=\id$ holds since
$m_\infty=m_1\ldots m_r$. For the same reason, any monodromy
$m_{j_0}$ can be expressed in terms of~$m_\infty$ and
the other monodromies~$m_j$, $j\ne j_0$; hence, the corresponding
relation can be omitted.
\endproof

\Remark
Note that, unlike, \eg,~\cite{degt.e6}, where the case of a
nonsingular fiber at infinity is considered, here
the relation
$m_\infty=\id$ does \emph{not} automatically follow from the
relation at infinity. Both relations are computed
in~\ref{inf.e8}--\ref{inf.e6} below.
\endRemark

\Remark
Usually, it is convenient to take for $\{\Gg_1,\ldots,\Gg_r\}$ a
\emph{geometric} basis for the group of the punctured plane
$\Delta^\circ$: each basis element is represented by the loop
composed of the counterclockwise boundary of a small disk
surrounding a puncture and a simple arc connecting this disk to
the reference point; all disks and arcs are assumed pairwise
disjoint except at the common reference point.
\endRemark

\subsection{The braid monodromy}\label{s.monodromy}
The braid monodromy of a trigonal curve~$\B$
can be computed using its
dessin (skeleton in the case of a maximal curve); below, we cite a
few relevant results of~\cite{degt.kplets}.

Recall that the braid group~$\BG3$ can be defined as
$\<\Gs_1,\Gs_2\,|\,\Gs_1\Gs_2\Gs_1=\Gs_2\Gs_1\Gs_2\>$;
it acts on the free group $\<\Ga_1,\Ga_2,\Ga_3\>$ \via
$$
\Gs_1\:(\Ga_1,\Ga_2,\Ga_3)\mapsto(\Ga_1\Ga_2\Ga_1^{-1},\Ga_1,\Ga_3),
\quad
\Gs_2\:(\Ga_1,\Ga_2,\Ga_3)\mapsto(\Ga_1,\Ga_2\Ga_3\Ga_2^{-1},\Ga_2).
$$
Introduce also the element $\Gs_3=\Gs_1\1\Gs_2\Gs_1$ and
consider the indices of $\Gs_1$, $\Gs_2$, $\Gs_3$ as residues
modulo~$3$, so that $\Gs_4=\Gs_1$ \etc.

The center of~$\BG3$ is generated by $(\Gs_1\Gs_2)^3$.
We denote by~$\bar\Gb$ the image of a braid $\Gb\in\BG3$ in the
\emph{reduced braid group}
$\BG3/(\Gs_1\Gs_2)^3\cong\CG2\mathbin*\CG3$. A braid~$\Gb$ is
uniquely recovered from~$\bar\Gb$ and its \emph{degree}
$\deg\Gb\in\BG3/[\BG3,\BG3]=\Z$. Recall that $\deg\Gs_i=1$.

\Remark
To be consistent with~\cite{degt.kplets}, we use the left action
of~$\BG3$ on the free group
$\<\Ga_1,\Ga_2,\Ga_3\>$. It appears, however, that
the right action is more suitable for the braid
monodromy, as it makes the map $\pi_1(\Delta^\circ)\to\BG3$ a
homomorphism rather than an anti-homomorphism. One can check that,
with one exception, all expressions involving braids are
symmetric modulo the central element
$(\Gs_1\Gs_2)^3=(\Gs_2\Gs_1)^3$. Hence, the only change needed to
pass to the right action is the definition of~$\Gs_3$: it should
be defined \via\ $\Gs_1\Gs_2\Gs_1\1$.
\endRemark

As in~\ref{s.strategy}, we fix a disk $\Delta\subset\Cp1$ and a
proper section~$s$ over~$\Delta$. All vertices, paths, \etc\. below
are assumed to belong to~$\Delta$.

For a trivalent \black-vertex~$v$ of the skeleton~$\Sk$ of~$\B$,
let $F_v$ be the fiber over~$v$ and let
$\pi_v=\pi_1(F_v\sminus(\B\cup E),s(v))$. A marking at~$v$,
see~\ref{s.marking}, gives rise to a natural ordering
$\{p_1,p_2,p_3\}$ of the three points of the intersection
$F_v\cap\B$ and, hence, to a \emph{canonical basis}
$\{\Ga_1,\Ga_2,\Ga_3\}$ for~$\pi_v$, see Figure~\ref{fig.basis},
which is well
defined up to simultaneous conjugation of the generators
by a power of $\Ga_1\Ga_2\Ga_3$, \ie, up to the action of the
central element $(\Gs_1\Gs_2)^3\in\BG3$.
(In the figure, $\Ga_i$ is a
small loop about~$p_i$, $i=1,2,3$.)

\midinsert
\centerline{\picture{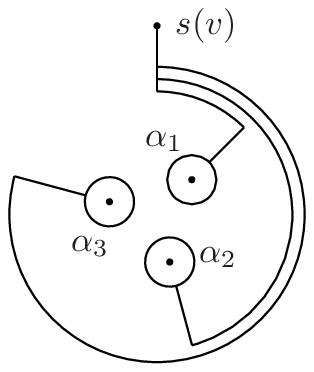}}
\figure
The canonical basis
\endfigure\label{fig.basis}
\endinsert

We extend the notion of canonical basis to the star of~$v$
\emph{in the dessin}, \ie, to the bold and solid edges incident
to~$v$,
extending
to but not including the \white-- and \cross-vertices. Over these
edges, the three points $\{p_1,p_2,p_3\}$ still form a proper
triangle, see~\cite{degt.kplets}; hence, they are still ordered by the
marking at~$v$ and one can construct the loops by combining radial
segments and arcs of a large circle. (Alternatively, one can
define this basis as the one obtained by translating a canonical
basis over~$v$ along the corresponding edge of the dessin.)

Given two marked trivalent \black-vertices~$u$ and~$v$
of~$\Sk$, one can
identify~$\pi_u$ and~$\pi_v$ by identifying the canonical bases
defined by the marking. This identification is well defined up to
the action of the center of~$\BG3$ (as so are the canonical
bases). If $u$ and~$v$ are connected by a path~$\Gg$, one can drag
the nonsingular fiber along~$\Gg$ and define the braid monodromy
$m_\Gg\:\pi_u\to\pi_v$. Combining this with the identification
above, one can define the element $\bar
m_\Gg\in\BG3/(\Gs_1\Gs_2)^3$. In particular, this construction
applies if $u$ and~$v$ are connected by an edge~$e$ of~$\Sk$;
depending on the type of the edge,
$\bar m_e$ is
given by the following expressions:
$$
\bar m_{[i,i+1]}=\bar\Gs_i,\quad \bar m_{[i+1,i]}=\bar\Gs_{i}^{-1},
 \quad\text{and}\quad
 \bar m_{[i,i]}=\bar\Gs_i\bar\Gs_{i-1}\bar\Gs_i.\eqtag\label{eq.bb}
$$
Using these relations, one can compute the
reduced monodromy~$\bar m_\Gg$
for any path~$\Gg$ composed of edges of~$\Sk$ connecting trivalent
\black-vertices. If $\Gg$ is a loop, the true monodromy~$m_\Gg$ is
recovered from~$\bar m_\Gg$ and the degree $\deg m_\Gg$, which
equals the total multiplicity of the singular fibers of~$\B$
encompassed by~$\Gg$. (The multiplicity of a
singular fiber~$F$ can be defined
as the number of the simplest type~$\tA{_0^*}$ fibers into which
$F$ can split.)

Now, let~$v$ be a marked
trivalent
\black-vertex of the dessin of~$\B$, and let~$u$ be
the \cross-vertex connected to~$v$
by the solid edge~$e'_i$.
Assume that the valency of~$u$ is~$d$, so that the singular
fiber~$F_u$ over~$u$ is of
type~$\tA{_{d-1}}$ ($\tA{_0^*}$ if $d=1$) or~$\tD{_{d+4}}$, see
Table~\ref{tab.j}.
Let~$\Gg$ be the loop composed of a small counterclockwise
circle around~$u$
connected to~$v$ along~$e'_i$.
Then, in any
canonical
basis for~$\pi_v$
defined by the marking at~$v$, the monodromy $m_\Gg$
along~$\Gg$ is given by
$$
\gathered
m_\Gg=\Gs_{i+1}^d,\
 \text{if $F_u$ is a type~$\tA{}$ fiber, or}\\
m_\Gg=\Gs_{i+1}^d(\Gs_1\Gs_2)^3,\
 \text{if $F_u$ is a type~$\tD{}$ fiber.}
\endgathered
\eqtag\label{eq.bc}
$$

\Remark
It is obvious geometrically (and can easily be shown formally)
that the reduced monodromy given by~\eqref{eq.bb} along the
boundary of a $d$-gonal region of the skeleton, when lifted to a
braid of appropriate degree, coincides with the monodromy about a
$d$-valent \cross-vertex given by~\eqref{eq.bc}.
\endRemark

\subsection{The relation at infinity and the monodromy at infinity}\label{s.infty}
We keep using the notation of~\ref{s.strategy}. In order to
compute the relation at infinity $[\partial\bGamma]=1$, assume
that $\Delta$ is the closure of the projection of the
disk~$\bGamma$ and that the reference fiber~$F_0$ is chosen
in~$\partial\Delta$.

Dragging a nonsingular fiber along~$\partial\bGamma$
and keeping
two points in~$s$ and~$\partial\bGamma$, one can define the
monodromy~$m$ along~$\partial\bGamma$ as an automorphism of the
relative homotopy set
$\pi_1(F_0\sminus(\B\cup E),(F_0\cap\partial\bGamma)\cup s(F_0),s(F_0))$.
Pick a
path~$p$ connecting $s(F_0)$ to $F_0\cap\partial\bGamma$ in
$F_0\sminus(\B\cup E)$. Then one has
$$
p\cdot[\partial\bGamma]\cdot m(p)\1=[s(\partial\Delta)]\1.
$$
This relation holds in any reasonable fundamental group, \eg, in
the group of the complement of $(\B\cup E)$ in the pull-back of
$\partial\Delta$.
Indeed, when dragged with the fiber, the arc~$p$ spans a
square~$S$
shown in Figure~\ref{fig.square}, disjoint from~$\B$ and~$E$,
and the product of the four
paths forming the boundary $\partial S$ (with appropriate
orientations) is a null homotopic loop. Note that the counterclockwise
directions of $\partial\bGamma$ and~$\partial\Delta$ (induced from
the complex orientations of the respective disks) are opposite to
each other.

\midinsert
\centerline{\picture{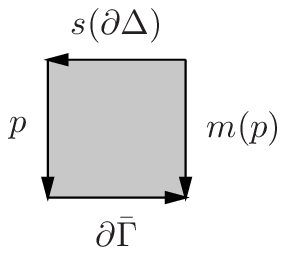}}
\figure\label{fig.square}
The square spanned by~$p$
\endfigure
\endinsert

Since $[s(\partial\Delta)]=1$ in
$\pi_1(\Sigma_k\sminus(\B\cup E\cup F))$ (the loop is
contractible along $s(\Delta)$),
the relation at infinity $[\partial\bGamma]=1$
takes the form
$$
p\cdot m(p)\1=1.
$$
The image $m(p)$ of a suitable arc~$p$ can easily be found using the
local forms given by~\eqref{eq.image.e8}, \eqref{eq.image.e7},
and~\eqref{eq.image.e6}. One can take for~$\Gamma$ a
small disk in
the line $v=d$, $d=\const$, so that $\bGamma$ is
the disk
$$
\{\ls|x|\le\epsilon,\ y=x^{k-1}\!/d\},
\eqtag\label{eq.Gamma}
$$
and compute the monodromy along the loop
$x=\epsilon\exp(2\pi it)$, $t\in[0,1]$. (Note that, in the
coordinates $(x,y)$, the `constant' section $s\:x'\mapsto c=\const$
is given by $x\mapsto cx^k$,
see~\eqref{eq.coord}; hence, the base point makes $k$ full
turns about the origin.) We omit the details, merely stating the
result in Subsections~\ref{inf.e8}--\ref{inf.e6} below.

One can use the same local models to compute the monodromy at
infinity~$m_\infty$. In other words, $m_\infty$ is the local
braid monodromy about~$F$ (in the clockwise direction) composed
with $(\Gs_1\Gs_2)^{3k}$ (due to the fact that the base point
makes $k$ full turns about the origin).
Below, we compute and simplify the
group of relations $[\partial\bGamma]=1$, $m_\infty=\id$, which
are present
in the fundamental group of any curve in question,
see Proposition~\ref{presentation}.

The results of the computation are stated in
Subsections~\ref{inf.e8}--\ref{inf.e6}.
We take for the reference fiber~$F_0$ the fiber $F_v$ over a
trivalent \black-vertex~$v$
of the dessin of~$\B$
connected by an edge to the vertex~$u$ corresponding
to the distinguished fiber~$F$,
and take for $\{\Ga_1,\Ga_2,\Ga_3\}$ a
canonical
basis
in~$\pi_F=\pi_v$ defined by a marking at~$v$. The particular
choice of the marking in each case is described below.

\subsection{The case of type~$\bE_8$}\label{inf.e8}
Assume that $P$ is of type~$\bE_8$, and hence $F$ is of type
$\tA{_1^*}$, see Proposition~\ref{1-1.e8}. Let~$\frak b$ be the
branch of~$\B$ at~$\PP$ that is not vertical.
Then $\bGamma$ is an ordinary
tangent to~$\frak b$, see~\eqref{eq.image.e8} and~\eqref{eq.Gamma},
and the relation at infinity is
$$
\Gr^3=\Ga_1\Ga_2^2,
\eqtag\label{eq.e8.1}
$$
assuming that $\Ga_2$ is represented by a loop about~$\frak b$
(so that the edge $[v,u]$ is~$e_2$ at~$v$).
The monodromy at infinity is
$m_\infty=(\Gs_1\Gs_2)^9(\Gs_1\Gs_2\Gs_1)\1$, and
the corresponding braid relations are
$$
\Ga_1=\Gr^{-3}(\Ga_1\Ga_2)\Ga_3(\Ga_1\Ga_2)\1\Gr^3,\quad
[\Ga_2,\Gr^{-3}\Ga_1]=1,\quad
\Ga_3=\Gr^{-3}\Ga_1\Gr^3.
$$
In view of~\eqref{eq.e8.1}, the second relations becomes
a tautology and the other two turn into
$$
\Ga_3=\Ga_2\Ga_1\Ga_2\1
\quad\text{and}\quad
[\Ga_1,\Ga_2^3]=1.
\eqtag\label{eq.e8.2}
$$
In particular, $\Ga_2^3$ is a central element.

\subsection{The case of type~$\bE_7$}\label{inf.e7}
Assume that $P$ is of type~$\bE_7$, and hence $F$ is of type
$\tA{_1}$, see Proposition~\ref{1-1.e8}. Then $\bGamma$ is
tangent to the distinguished branch~$\frak b$ of~$\B$ at~$\PP$.
Unless
stated otherwise, we will choose the basis $\Ga_1$, $\Ga_2$,
$\Ga_3$ so that
\roster
\item"$(*)$"
$\Ga_2$ and $\Ga_3$ are represented by
loops about the two branches of~$\B$ at~$\PP$ and
$\Ga_2$
corresponds to the distinguished branch~$\frak b$.
In particular, $[v,u]$ is the edge~$e'_1$ at~$v$.
\endroster
(Occasionally, we will also consider the case when
the generator corresponding to~$\frak b$ is
$\Ga_3$.)
Then,
the relation at infinity is
$$
\Gr^3=\Ga_2\Ga_3\Ga_2
\quad\text{or}\quad
\Gr^3=\Ga_2\Ga_3\Ga_2\Ga_3\Ga_2\1,
$$
assuming that $\Ga_2$
or, respectively, $\Ga_3$
corresponds to~$\frak b$. The monodromy at infinity
is $m_\infty=(\Gs_1\Gs_2)^9\Gs_2^{-2}$,
the corresponding braid relations being
$$
[\Ga_1,\Gr^3]=1
\quad\text{and}\quad
[\Ga_i,\Gr^{-3}(\Ga_2\Ga_3)]=1,\ i=2,3.
$$
Combining the last pair of relations
with the relation at infinity, one concludes that
(assuming that $\Ga_2$ corresponds to~$\frak b$)
$$
[\Ga_2,\Ga_3]=1
\quad\text{and}\quad
[\Ga_i,\Gr^3]=[\Ga_i,\Ga_2^2\Ga_3]=1,
\ i=1,2,3;
\eqtag\label{eq.e7.1}
$$
then, the relation at infinity takes the form
$$
\Gr^2\Ga_1=\Ga_2.
\eqtag\label{eq.e7.2}
$$
If the generator corresponding to~$\frak b$
is~$\Ga_3$, instead of \eqref{eq.e7.1} and~\eqref{eq.e7.2}
one has
$$
[\Ga_2,\Ga_3]=[\Ga_i,\Gr^3]=[\Ga_i,\Ga_2\Ga_3^2]=1,\ i=1,2,3,
\quad\text{and}\quad
\Gr^2\Ga_1=\Ga_3.
\eqtag\label{eq.e7.alt}
$$

\subsection{The case of type~$\bE_6$}\label{inf.e6}
Assume that $P$ is of type~$\bE_6$, and hence $\PP$ is of type
$\bA_5$, see Proposition~\ref{1-1.e8}. Then $\bGamma$ is
inflection tangent to each of the two branches of~$\B$ at~$\PP$,
and the relation at infinity is
$$
\Gr^4=(\Ga_2\Ga_3)^3,
\eqtag\label{eq.e6}
$$
assuming that $\Ga_2$ and $\Ga_3$
are represented by
loops about the two branches at~$\PP$
(so that $[v,u]$ is the edge~$e'_1$ at~$v$).
The monodromy at infinity is
$m_\infty=(\Gs_1\Gs_2)^{12}\Gs_2^{-6}$,
and the corresponding braid relations are
$$
[\Ga_1,\Gr^4]=1
\quad\text{and}\quad
[\Ga_i,\Gr^{4}(\Ga_2\Ga_3)^{-3}]=1,\ i=2,3.
$$
These relations follow from~\eqref{eq.e6}.

\section{The inclusion homomorphism}\label{S.inclusion}

Here, we compute the homomorphism of the fundamental
groups induced by the inclusion to~$\Cp2$ of a Milnor ball~$\MB$
about a
type~$\bE$ singular point~$P$ of a sextic~$B$. These results are used
in~\S\ref{S.perturbations}.

\subsection{The setup}
In order to compute the inclusion homomorphisms, we represent the
sextic~$B$ by the polynomial given by~\eqref{eq.curve.e8},
\eqref{eq.curve.e7}, or~\eqref{eq.curve.e6}, assuming all
parameters involved real and positive, and generate the group
$\pi_1(\MB\sminus B)$ by the classes of appropriately chosen
loops $\bc_1$, $\bc_2,\bc_3$ in
the complement $\{v=\epsilon\}\sminus B$, where $\epsilon\ll1$ is
a positive real constant. (In what follows, we identify the loops
and their classes.) Each loop~$\bc_i$, $i=1,2,3$, is composed of a
small circle~$C_i$ about a point of intersection
$\{v=\epsilon\}\cap B$ and a path~$p_i$ connecting
a point $r_i\in C_i$ to
the base point, which is a large real number.

The image of the line $\{v=\epsilon\}$ in~$\Sigma_k$ is the
parabola $\{x^{k-1}=\epsilon y\}$. It intersects the `constant'
section $\{y'=c\}=\{y=cx^k\}$ used to define the braid monodromy
at the origin and at the point
$r_0:=(x_0,y_0)=(1/\epsilon c,1/\epsilon^k c^{k-1})$.
We assume that~$c$ is also
a real constant, $0\ll c\ll1/\epsilon$, so that $y_0\gg0$,
and take~$r_0$
for the common base point in both
the line
$\{v=\epsilon\}$ and the reference fiber
$F'_0=\{x=x_0\}$. (This fiber may differ from the
reference fiber
considered in~\S\ref{S.group};
the necessary adjustments are explained below.)

Now, consider the fiber~$F'_i$, $i=1,2,3$, passing through~$r_i$.
(We keep the same notation~$C_i$, $r_i$, and~$p_i$ for the images
of the corresponding elements in~$\Sigma_k$.) The point~$r_i$ is
close to a branch of~$\B$; let
$\bc'_i\in\pi_1(F'_i\sminus(\B\cup E),r_i)$ be the element
represented by a small circle through~$r_i$ encompassing this
branch. Dragging~$F'_i$ along~$p_i$ while keeping the base point
in~$p_i$, one defines the braid monodromy
$$
m'_i\:\pi_1(F'_i\sminus(\B\cup E),r_i)\to
 \pi_1(F'_0\sminus(\B\cup E),r_0).
$$
(One should make sure that $p_i$ does not pass through the origin
in the line $\{v=\epsilon\}$.)
It is immediate that $m'_i(\bc'_i)$ represents the image of the
generator~$\bc_i$ under the inclusion homomorphism (\cf.
Figure~\ref{fig.incl}, where the curve $\B$, the line
$\{v=\epsilon\}$, and the section~$s$ are drawn in bold, dashed,
and dotted lines, respectively; the two grey lassoes, one in
$\{v=\epsilon\}$ and one in the fiber $F_i'=F_0'$,
represent the
same element of the fundamental group).

\midinsert
\centerline{\picture{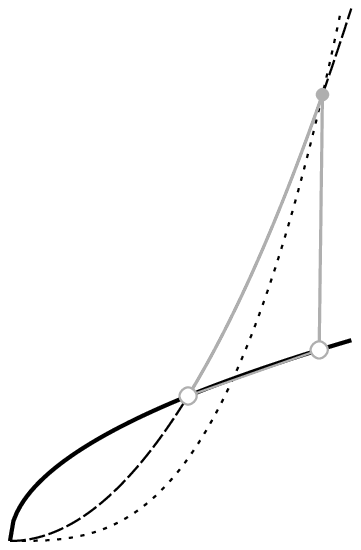}}
\figure\label{fig.incl}
Computing the inclusion homomorphism
\endfigure
\endinsert

\subsection{The case of type~$\bE_7$}\label{incl.e7}
The original curve~$B$ is given by~\eqref{eq.curve.e7}. All
three points of intersection of~$B$ and the line $\{v=\epsilon\}$
are real, and we take for $\{\bc_1,\bc_2,\bc_3\}$ a `linear'
basis, numbering the intersection points consecutively by the
decreasing of the $u$-coordinates and taking for~$p_i$ segments of
the real line, circumventing the interfering intersection points
and the origin
in the counterclockwise direction.

All three points of intersection of~$\B$ and the reference
fiber~$F_0$ are also real, and we choose a similar `linear' basis
$\{\Ga'_1,\Ga'_2,\Ga'_3\}$ for the group
$\pi_1(F'_0\sminus(\B\cup E),r_0)$. Then one has
$$
\bc_1\mapsto\Ga'_1,\quad
\bc_2\mapsto\Ga'_2,\quad
\bc_3\mapsto(\Ga'_2\Ga'_3)\1\Ga'_1(\Ga'_2\Ga'_3).
$$
In order to pass to the reference fiber~$F_0$ considered
in~\ref{inf.e7} and a
canonical
basis $\{\Ga_1,\Ga_2,\Ga_3\}$
satisfying~\ref{inf.e7}$(*)$, one can drag~$F'_0$ along the arc
$x=x_0\exp(\pi it/2)$, $t\in[0,1]$. Then $\Ga'_1=\Ga_1$,
$\Ga'_2=\Ga_2\Ga_3\Ga_2\1=\Ga_3$ (we use the commutativity
relation in~\eqref{eq.e7.1}), and $\Ga'_3=\Ga_2$. Finally, the
inclusion homomorphism is given by
$$
\bc_1\mapsto\Ga_1,\quad
\bc_2\mapsto\Ga_3,\quad
\bc_3\mapsto(\Ga_2\Ga_3)\1\Ga_1(\Ga_2\Ga_3).
\eqtag\label{eq.incl.e7}
$$

\subsection{The case of type~$\bE_8$}\label{incl.e8}
The curve~$B$ is given by~\eqref{eq.curve.e8}, and the
points of intersection of~$B$ and the line $\{v=\epsilon\}$ are
the three roots $u=\root3\of{\epsilon^5}$. Let
$\{\bc_1,\bc_2,\bc_3\}$ be a basis similar to the one shown in
Figure~\ref{fig.basis}, with the paths $p_i$ composed of radial
segments and the arcs $u=\const\cdot\exp(\pm2\pi it/3)$,
$t\in[0,1]$. We assume that the generator corresponding to the
real branch of~$B$ is~$\bc_2$.

All three points of intersection of~$\B$ and the reference
fiber~$F_0$ are real, and we choose a `linear' basis
$\{\Ga_1,\Ga_2,\Ga_3\}$ as in~\ref{incl.e7} for the group
$\pi_1(F'_0\sminus(\B\cup E),r_0)$. Dragging~$F_0$ along the arc
$x=x_0\exp(\pi it)$, $t\in[0,1]$, to the fiber~$\{x=-\epsilon\}$,
one can see that
$\Ga_1$, $\Ga_2$, $\Ga_3$ are indeed equal to the basis
elements considered in~\ref{inf.e8}.

In these bases, the inclusion homomorphism is given by
$$
\bc_1\mapsto(\Ga_1\Ga_2)\Ga_3(\Ga_1\Ga_2)\1,\quad
\bc_2\mapsto\Ga_1,\quad
\bc_3\mapsto\Ga_3.
\eqtag\label{eq.incl.e8}
$$

\subsection{The case of type~$\bE_6$}\label{incl.e6}
The curve~$B$ is given by~\eqref{eq.curve.e6},
the intersection points of~$B$ and $\{v=\epsilon\}$ are
the three roots $u=\root3\of{\epsilon^5}$, and we take for
$\{\bc_1,\bc_2,\bc_3\}$ the same basis as in~\ref{incl.e8}.
The basis $\{\Ga_1,\Ga_2,\Ga_3\}$ in the reference fiber is chosen
`linear' as above; these elements do satisfy the conditions imposed
in~\ref{inf.e6}. In these bases, the inclusion homomorphism is
given by
$$
\bc_1\mapsto(\Ga_1\Ga_2\Ga_3)\Ga_1(\Ga_1\Ga_2\Ga_3)\1,\quad
\bc_2\mapsto\Ga_1,\quad
\bc_3\mapsto(\Ga_2\Ga_3)\1\Ga_1(\Ga_2\Ga_3).
\eqtag\label{eq.incl.e6}
$$

\section{The computation}\label{S.computation}

In this section, Theorems~\ref{th.e7} and~\ref{th.e7.group} are
proved. Throughout the section, we fix the following notation:
$B$ stands for an irreducible maximal plane sextic
with a distinguished type~$\bE_7$ singular point~$P$, and
$L\subset\Cp2$ is the line tangent to~$B$ at~$P$.
We denote by~$\B$ and~$F$,
respectively, the trigonal curve corresponding to~$B$ and its
distinguished fiber, see Proposition~\ref{1-1.e7}; $\Sk$ stands
for the skeleton of~$\B$.

\subsection{Proof of Theorem~\ref{th.e7}}\label{s.proof}
According to Proposition~\ref{sing.points}, all triple points
of~$B$ other than~$P$ are of type~$\bE_6$ or~$\bE_8$. Let $t=0$,
$1$, or~$2$ be their number. Then $\Sk$ has the following
properties:
\roster
\item\local{Sk.1}
$\deg\Sk=9-3t$ and
$\Sk$ has exactly $t$ singular vertices, none of which is \white-;
\item\local{Sk.2}
if $t=0$, then $\Sk$ does not admit a splitting marking
(Proposition~\ref{splitting}).
\endroster
Conversely, in view of Theorem~\ref{th.trigonal},
any skeleton~$\Sk$ satisfying~\loccit{Sk.1}
and~\loccit{Sk.2}
above (for some
integer $t\ge0$) represents an irreducible
maximal trigonal curve~$\B$ as in
Proposition~\ref{1-1.e7}; hence, it represents two
irreducible maximal sextics with a distinguished
type~$\bE_7$ singular point.

The distinguished fiber~$F$ is located at the center of a bigonal
region~$R$ of~$\Sk$. In the drawings below, we show the boundary
of~$R$ in grey.

Assume that all \black-vertices in the boundary of~$R$ are
nonsingular.
Then~$R$ looks as shown in Figure~\ref{fig.insertion}.
This fragment of the skeleton (if present)
is called the \emph{insertion}. The two branches of~$\B$ at the
node located in~$F$ are in a natural correspondence with the two
edges of~$\partial R$; hence, selecting one of the branches (see
Proposition~\ref{1-1.e7}) can be interpreted geometrically as
selecting one of the two arcs in the boundary of the insertion.

\midinsert
\centerline{\picture{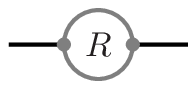}}
\figure\label{fig.insertion}
\endfigure
\endinsert

Removing the insertion and patching it with an edge, one obtains
another valid skeleton $\Sk'$ of degree $6-3t\le6$, \cf\.
Figure~\ref{fig.e7}. (Note that one has $t\le1$ in this case.)
The
\black-vertices
of~$\Sk'$ are in a one-to-one correspondence with the same valency
vertices of~$\Sk$ other than the two vertices in~$\partial R$.
Conversely, given a skeleton $\Sk'$ of degree $6-3t\le6$ with $t$
singular \black-vertices, $t=0,1$, one can place the insertion at the
middle of any edge of~$\Sk'$ to obtain a new skeleton $\Sk$
satisfying~\loccit{Sk.1} above.

\lemma
If $t=0$, the skeleton~$\Sk$ admits a splitting marking if and
only if so does~$\Sk'$.
\endlemma

\proof
It is
immediate that any splitting marking of~$\Sk$ restricts
to a splitting marking of~$\Sk'$ and, \viceversa, any splitting
marking of~$\Sk'$ extends (uniquely) to a splitting marking
of~$\Sk$.
\endproof

\Remark
This trick, replacing a given skeleton~$\Sk$ by
another skeleton~$\Sk'$ obtained from~$\Sk$ by removing a certain
fragment, appears on numerous occasions in the classification of
plane sextics and, more generally, in the study of extremal
elliptic surfaces. It would be interesting to understand if the
passage from~$\Sk$ to~$\Sk'$ corresponds to a simple geometric
construction defined in terms of trigonal curves or covering
elliptic surfaces. At present, I do not know any geometric
interpretation.
\endRemark

Thus, the classification of skeletons~$\Sk$ satisfying
conditions~\loccit{Sk.1}
and~\loccit{Sk.2} above
and containing an insertion can be done in two steps:
\Dashes
\dash
The classification of skeletons of degree~$6$ without singular
vertices and not admitting a splitting marking, and the
classification of skeletons of degree~$3$ with exactly one
singular vertex, which is~\black-. This is done
in~\cite{symmetric}, and the complete list is presented in
Figure~\ref{fig.e7}(a)--(e).
\dash
Placing an insertion with one of the two arcs selected
to one of the
edges of each skeleton~$\Sk'$ discovered at step one.
\endDashes
The second step is clearly equivalent to choosing a pair
$(e,\goth o)$, where $e$ is an edge of~$\Sk'$ and $\goth o$ is a
coorientation of~$e$. Such pairs are to be considered up to
orientation preserving automorphisms of~$\Sk'\subset S^2$.
All essentially distinct edges~$e$ (with the
coorientation~$\frak o$ ignored) are shown in
Figures~\ref{fig.e7}(a)--(e). Taking into account the coorientation,
one obtains the list given by
Table~\ref{tab.e7}, lines~$1$--$9$.

\midinsert
\centerline{\vbox{\halign{\hss#\hss&&\qquad\hss#\hss\cr
\cpic{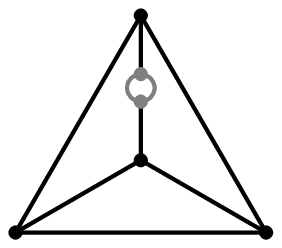}&
\cpic{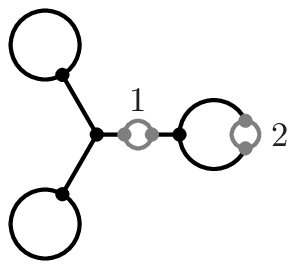}&
\cpic{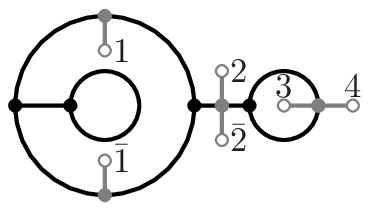}\cr
\noalign{\medskip}
(a)&(b)&(c)\cr
\crcr}}}
\bigskip
\centerline{\vbox{\halign{\hss#\hss&&\kern3em \hss#\hss\cr
\cpic{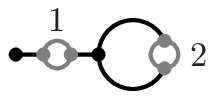}&
\cpic{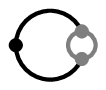}&
\cpic{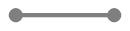}&
\cpic{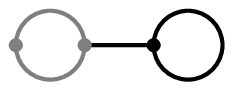}\cr
\noalign{\medskip}
(d)&(e)&(f)&(g)\cr
\crcr}}}
\figure\label{fig.e7}
\endfigure
\endinsert

The few remaining cases, when the boundary of~$R$ contains
singular \black-vertices, can easily be treated manually using the
vertex count given by Corollary~\ref{vertex.count}. The two
skeletons obtained are shown in Figure~\ref{fig.e7}(f), (g), and
the corresponding sets of singularities are listed in
Table~\ref{tab.e7}, lines~$10$, $11$.

\midinsert
\table\label{tab.e7}
Maximal sets of singularities with a type~$\bE_7$ point
represented by irreducible sextics
\endtable

\def\no{}
\def\*{\llap{$^*$\,}}
\def\tref#1{\,\text{\ref{#1}}\,}
\centerline{\vbox{\offinterlineskip\halign{%
\tabstrut\ \hss#\hss\ \vrule
&\quad$#$\hss\quad\vrule
&\quad\null#\hss\quad\vrule
&&\ \hss$#$\hss\ \vrule\cr
\noalign{\hrule}
\exstrut&&&&&\cr
\#&\text{Set of singularities}&
 \hss Figure&\text{Count}&
 \pi_1&S^\perp\cr
\exstrut&&&&&\cr
\noalign{\hrule}
\exstrut&&&&&\cr
1&\bE_7\splus2\bA_4\splus2\bA_2&
 \ref{fig.e7}(a)&
 (1,0)&
 \tref{s.no.1}&(15,0,15)\cr
2&\*\bE_7\splus\bA_{12}&
 \fragment(b)1&
 (0,1)&
 \tref{s.+loop}&(7,2,2)\cr
3&\*\bE_7\splus\bA_{10}\splus\bA_2&
 \fragment(b)2&
 (2,0)&
 \tref{s.loop}&(11,0,3)\cr
4&\bE_7\splus2\bA_6&
 \fragment(c){1,\bar1\!\!\!}&
 (0,1)&
 \tref{s.no.4}&(7,0,7)\cr
5&\*\bE_7\splus\bA_8\splus\bA_4&
 \fragment(c)2&
 (0,1)&
 \tref{s.+loop}&(23,2,2)\cr
6&\*\bE_7\splus\bA_6\splus\bA_4\splus\bA_2&
 \fragment(c)3&
 (2,0)&
 \tref{s.loop}&(35,0,3)\cr
7&\*\bE_7\splus\bE_6\splus\bA_6&
 \fragment(d)1&
 (0,1)&
 \tref{s.+loop}&(11,2,2)\cr
8&\*\bE_7\splus\bE_6\splus\bA_4\splus\bA_2&
 \fragment(d)2&
 (2,0)&
 \tref{s.loop}&(15,0,3)\cr
9&\bE_7\splus\bE_8\splus2\bA_2&
 \ref{fig.e7}(e)&
 (1,0)&
 \tref{s.loop}&(3,0,3)\cr
10&\bE_7\splus2\bE_6&
 \ref{fig.e7}(f)&
 (1,0)&
 \tref{s.no.10}&(3,0,3)\cr
11&\*\bE_7\splus\bE_8\splus\bA_4&
 \ref{fig.e7}(g)&
 (0,1)&
 \tref{s.+loop}&(3,2,2)\cr
\exstrut&&&&&\cr
\noalign{\hrule}
\crcr}}}
\endinsert

The results of the computation are collected in
Table~\ref{tab.e7}, where we list the set of singularities, the
skeleton~$\Sk$ of~$\B$, the number of deformation classes (see
below), and a reference to the computation of the fundamental
group.
A set of singularities is marked with a~${}^*$ if it is realized
by two equisingular deformation classes which have the same
skeleton but differ by the selected branch of the insertion. (In
the terminology of Proposition~\ref{1-1.e7}, the two families
differ by the distinguished branch of~$\B$ at~$\PP$.)
For completeness, we also list the lattice~$S^\perp$ corresponding
to the homological type of the sextic,
see~\cite{JAG} for the definitions:
the notation $(a,b,c)$ stands for the quadratic
form generated by two elements $u$, $v$ with
$u^2=2a$, $u\cdot v=b$, and $v^2=2c$. The lattice is obtained by
comparing two independent classifications, those of curves and of
abstract homological types, and taking into account the number of
classes obtained (see also Shimada~\cite{Shimada}).

The number of deformation classes is listed in the form
$(n_r,n_c)$, where $n_r$ is the number of real curves and $n_c$ is
the number of pairs of complex conjugate curves. (Thus, the total
number of classes is $n_r+2n_c$.) Real are the curves whose
skeletons admit an orientation reversing automorphism of
order~$2$ (with the marked arc taken into account); otherwise, two
symmetric skeletons represent a pair of complex conjugate curves.
\qed

\subsection{Proof of Theorem~\ref{th.e7.group}}\label{proof.e7.group}
We compute the fundamental groups using the strategy outlined
in~\ref{s.strategy} and the presentation given by
Proposition~\ref{presentation}.
As in~\ref{s.infty}, we choose for the reference fiber~$F_0$ the
fiber~$F_v$ over a \black-vertex~$v$ connected by an edge to the
\cross-vertex corresponding to~$F$, and consider a
canonical
basis
$\{\Ga_1,\Ga_2,\Ga_3\}$ for the group $\pi_F=\pi_v$ defined by an
appropriate marking. In most cases, we assume that the basis
satisfies~\ref{inf.e7}$(*)$. Then,
for all groups, the relations
$m_\infty=\id$ and $[\partial\bGamma]=1$ are given
by~\eqref{eq.e7.1} and~\eqref{eq.e7.2}, and the remaining braid
relations $m_j=\id$ are computed using the techniques outlined
in~\ref{s.monodromy}; in most cases, just a few extra
relations suffice
to show that the group is abelian. A detailed case by case
analysis is given in \ref{s.+loop}--\ref{s.no.1} below.

A great deal of the calculation in~\ref{s.+loop}--\ref{s.no.1} was
handled using \GAP~\cite{GAP}. In most cases, we merely input the
relations and query the size of the resulting group; having
obtained
six, we know that the group is~$\CG6$. In the more advanced case
in~\ref{s.no.1}, we quote the \GAP\ input/output in
Diagram~\ref{GAP.no.1}.

\subsection{Sets of singularities Nos.~2, 5, 7, and 11}\label{s.+loop}
Assume that the distinguished fiber $F$ has a neighborhood shown
in Figure~\ref{fig.e7-stem}. (The leftmost \black-vertex can be
either bi- or trivalent; it is not used in the calculation.)

\midinsert
\centerline{\picture{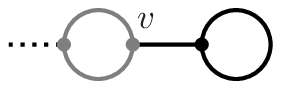}}
\figure\label{fig.e7-stem}
\endfigure
\endinsert

Assume that
the distinguished branch is such that
the basis in~$\pi_v$
satisfies~\ref{inf.e7}$(*)$.
(The case when the generator
corresponding to the distinguished branch is $\Ga_3$
is treated similarly;
alternatively, one can argue that the corresponding fragment
is obtained from the one considered by the complex
conjugation.) Then, the braid relation about the type~$\tA{_0^*}$
singular fiber represented by the rightmost loop of the skeleton
is $\Ga_2=(\Ga_2\1\Ga_1\Ga_2)\Ga_3(\Ga_2\1\Ga_1\Ga_2)\1$.
In the presence of~\eqref{eq.e7.1},
it simplifies to
$\Ga_1\1\Ga_2\Ga_1=\Ga_3$.
Let
$$
G=\<\Ga_1,\Ga_2,\Ga_3\,|\,
 \text{\eqref{eq.e7.1}, \eqref{eq.e7.2}, $\Ga_1\1\Ga_2\Ga_1=\Ga_3$}\>.
\eqtag\label{eq.+loop}
$$
Since $\Ga_2$, $\Ga_3$ commute and
$\Ga_1\1(\Ga_2^2\Ga_3)\Ga_1=(\Ga_2^2\Ga_3)$, one has
$\Ga_1\1\Ga_3\Ga_1=\Ga_2^2\Ga_3\1$. Thus, the conjugation
by~$\Ga_1$ preserves the abelian subgroup generated by~$\Ga_2$
and~$\Ga_3$, and the map $t\:w\mapsto\Ga_1\1w\Ga_1$ is given by
$$
t\:\Ga_2\mapsto\Ga_3\mapsto\Ga_2^2\Ga_3\1\mapsto\Ga_3^3\Ga_2^{-2}
 \mapsto\ldots.
$$

Using the non-commutativity relations obtained, one can move
all three copies of~$\Ga_1$ in~\eqref{eq.e7.2} to the left; this
gives $\Ga_1^3=\Ga_2\1\Ga_3^{-2}$. In particular, $t^3=\id$ and
hence $\Ga_2^3=\Ga_3^3$. Now, it is easy to see that the
commutant $[G,G]\cong\CG3$ is generated by the central order~$3$
element $\Ga_2\Ga_3\1$ and the abelianization
$\bar G:=G/[G,G]$ is the group
$$
\operatorname{ab}\<\bGa_1,\bGa_2,\bGa_3\,|\,
 \bGa_2=\bGa_3,\ 3(\bGa_1+\bGa_2)=0\>.
$$
(Here, the bar stands for the projection of an element to the
abelianization.)

For any group~$\pi_1$ with central commutant (for example, for
any quotient group of~$G$ above), the map
$\bar x\wedge\bar y\mapsto[x,y]\in[\pi_1,\pi_1]$ is a well defined
skew-symmetric bilinear form
$\bigwedge^2\bar\pi_1\to[\pi_1,\pi_1]$ (where
$\bar\pi_1:=\pi_1/[\pi_1,\pi_1]$); clearly, the group is abelian
if and only if this form is identically zero.
In particular, abelian are the
fundamental groups of irreducible sextics with the sets of
singularities Nos.~2, 5, 7, and 11 in Table~\ref{tab.e7}, as in
this case $\bar\pi_1$ is cyclic and $\bigwedge^2\bar\pi_1=0$.

\subsection{Proof of Proposition~\ref{th.split} and Corollary~\ref{th.split.pert}}\label{proof.split}
The arguments in~\ref{s.+loop} apply as well to
a \emph{reducible} maximal sextic $B$,
provided
that the skeleton of the trigonal model of~$B$
has a fragment shown
in Figure~\ref{fig.e7-stem}. Such skeletons can be enumerated
similar to~\ref{s.proof}, placing an insertion to an appropriate
edge of one of the skeletons of reducible trigonal curves found
in~\cite{symmetric}. There are four skeletons with a loop, which
result in the five sets of singularities listed in
Proposition~\ref{th.split}. (In one case, there is an extra choice
of the singular fiber to be converted to a $\bD$-type point
by an elementary transformation.)

Figure~\ref{fig.e7-stem} implies that the two branches of~$\B$
at~$\PP$ are in the same irreducible component; hence, the
corresponding sextic~$B$ splits into two irreducible cubics (at
least one having a cusp). Furthermore, analyzing the other
possibilities, one can deduce that each set of singularities
listed in Proposition~\ref{th.split} is realized by a unique, up
to complex conjugation, equisingular deformation family of sextics
\emph{splitting into two cubics}. (The set of singularities
$\bE_7\splus\bA_9\splus\bA_2\splus\bA_1$ is also realized by two
other families:
one can start from the skeleton marked
$\tilde\bA_7\splus\tilde\bA_1\splus2\tilde\bA_0^*$
in~\cite{symmetric} and place
the insertion as shown in Figure~\ref{fig.e7-leaf}.
In this case, the two branches at~$\PP$ are
in two distinct components of~$\B$ and,
depending on the branch distinguished,
the corresponding sextic splits into a quintic and a
line or a quartic and a conic.)

As explained in~\ref{s.+loop},
the fundamental group~$\pi_1$ of each of the curves obtained
is a quotient of the group~$G$ given
by~\eqref{eq.+loop},
and the commutator form $\bigwedge^2\bar\pi_1\to[\pi_1,\pi_1]$ is
determined by the value $[\bGa_1,\bGa_2]$.
Thus, the only question
is whether an extra relation in the fundamental group implies
$[\bGa_1,\bGa_2]=0$; if such a relation does exist,
$\pi_1$
is abelian; otherwise,
$\pi_1=G$.
The remaining relations are easily found
using the techniques
explained in~\ref{s.monodromy}.
For example, the monodromy along the boundary of
the $2n$-gon adjacent to the insertion (assuming that this region
does not contain a $\tD{}$-type fiber)
gives the relation
$(\Ga_1\Ga_2)^n=(\Ga_2\Ga_1)^n$; in~$G$, it simplifies to
$n[\bGa_1,\bGa_2]=0$; hence, it is a tautology if $n=0\bmod3$, and
it implies that the quotient is abelian if $n\ne0\bmod3$.
(Note that this region must have an even number of corners, as
otherwise the relation would imply $\bGa_1=\bGa_2$ and the curve
would be irreducible.)
We omit further details; the final result is stated in
Proposition~\ref{th.split}.
\qed

\subsection{Sets of singularities Nos.~3, 6, 8, and 9}\label{s.loop}
Assume that the distinguished fiber $F$ has a neighborhood shown
in Figure~\ref{fig.e7-leaf}. (The leftmost \black-vertex can be
either bi- or trivalent; it is not used in the calculation.) In
other words, we assume that one of the regions adjacent to the
bigon~$R$ containing~$F$ is a triangle. Then, from Figure~\ref{fig.e7}
it follows that the other region adjacent to~$R$ is a $3$-,
$5$-, $7$-, or $11$-gon.

\midinsert
\centerline{\picture{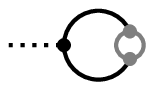}}
\figure\label{fig.e7-leaf}
\endfigure
\endinsert

Over one of the two grey vertices in Figure~\ref{fig.e7-leaf}
(depending on the distinguished branch) a
canonical
basis can be chosen to satisfy~\ref{inf.e7}$(*)$.
In this basis,
in addition to~\eqref{eq.e7.1} and~\eqref{eq.e7.2}, the singular
fibers inside the two regions adjacent to~$R$ give the relations
$$
(\Ga_1\Ga_2)^m\Ga_1=\Ga_2(\Ga_1\Ga_2)^m,
\quad
(\Ga_1\Ga_3)^n\Ga_1=\Ga_3(\Ga_1\Ga_3)^n,
\eqtag\label{eq.leaf}
$$
where either $m=1$ and $n=1,2,3,5$ or $n=1$ and $m=1,2,3,5$.
(The relations are given by~\eqref{eq.bc};
for the second relation in~\eqref{eq.leaf},
we use the commutativity $[\Ga_2,\Ga_3]=1$.) Using
\GAP~\cite{GAP}, one can see that, for each pair $(m,n)$ as above,
the group
$\<\Ga_1,\Ga_2,\Ga_3\,|\,\text{\eqref{eq.e7.1}, \eqref{eq.e7.2},
 \eqref{eq.leaf}}\>$
has order six; hence, it is abelian.

\subsection{The set of singularities $\bE_7\splus2\bE_6$ (No.~10)}\label{s.no.10}
In this special case, the skeleton of~$\B$
has no trivalent \black-vertices.
However, we can choose a basis similar to the one shown in
Figure~\ref{fig.basis} and satisfying~\ref{inf.e7}$(*)$ in the
fiber~$F_0$ over a point in the (open) solid edge connecting the
\cross-vertex representing~$F$ and an appropriate \black-vertex.
Then, both~\eqref{eq.e7.1} and~\eqref{eq.e7.2} still hold and, in
addition,
the monodromy about the
type~$\tE{_6}$ singular fiber close to~$F_0$ gives the relations
$$
\Ga_2\Ga_1\Ga_2\Ga_3=\Ga_1\Ga_2\Ga_3\Ga_1,
\quad
\Ga_3\Ga_1\Ga_2\Ga_3=\Ga_1\Ga_2\Ga_3\Ga_2.
$$
Using \GAP~\cite{GAP}, one can conclude that the resulting group has
order six; hence, it is abelian.

\subsection{The set of singularities $\bE_7\splus2\bA_6$
(No.~4)}\label{s.no.4}
After an automorphism and/or complex conjugation, one can assume
that the insertion~$R$
is as shown in Figure~\ref{fig.2a4} and
that a
canonical
basis over~$v$ can be chosen to
satisfy~\ref{inf.e7}$(*)$.

\midinsert
\centerline{\picture{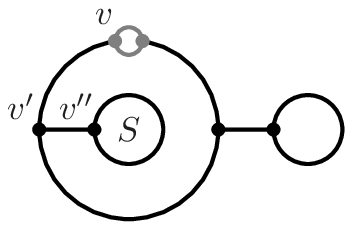}}
\figure\label{fig.2a4}
The set of singularities $\bE_7\splus2\bA_6$
\endfigure
\endinsert

Then,
in addition to~\eqref{eq.e7.1} and~\eqref{eq.e7.2}, one has
relations~\eqref{eq.leaf} with $m=n=3$ (given by the two heptagons
adjacent to the insertion)
and an additional relation $\Ga_1'=\Ga_2'\Ga_3'{\Ga_2'}\1$
given by the monodromy along
the loop $[\partial S]$, where~$S$ is the region shown
in Figure~\ref{fig.2a4}
and
$\Ga_1'$, $\Ga_2'$, $\Ga_3'$ is an appropriate
canonical
basis in the fiber over~$v''$. Connecting~$v''$ to~$v$ \via~$v'$
and using~\eqref{eq.bb}, one
obtains
$$
\Ga_1'=(\Ga_1\Ga_2)\1\Ga_1(\Ga_1\Ga_2),
\quad
\Ga_2'=(\Ga_1\Ga_2)\1\Ga_2(\Ga_1\Ga_2),
\quad
\Ga_3'=\Ga_3,
$$
and, in view of the fact that $\Ga_2$ and~$\Ga_3$ commute, the
relation
about~$[\partial S]$
simplifies to
$$
\Ga_1=(\Ga_2\Ga_1)\Ga_3(\Ga_2\Ga_1)\1.
$$
Using \GAP~\cite{GAP}, one finds that the group obtained has
order six; hence, it is abelian.

\subsection{The set of singularities $\bE_7\splus2\bA_4\splus2\bA_2$ (No.~1)}\label{s.no.1}
After an automorphism and/or complex conjugation, one can assume
that the skeleton of~$\B$
is as shown in Figure~\ref{fig.4a2} and
that a
canonical
basis over~$v$ can be chosen to
satisfy~\ref{inf.e7}$(*)$.

\midinsert
\centerline{\picture{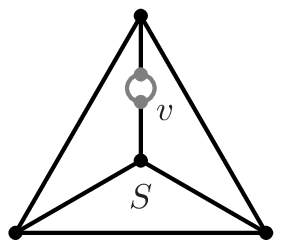}}
\figure\label{fig.4a2}
The set of singularities $\bE_7\splus2\bA_4\splus2\bA_2$
\endfigure
\endinsert

Then,
in addition to~\eqref{eq.e7.1} and~\eqref{eq.e7.2}, one has
relations~\eqref{eq.leaf} with $m=n=2$ (given by the two pentagons
adjacent to the insertion)
and the relation
$$
(\Ga_1\Ga_3\Ga_1\1)\Ga_2(\Ga_1\Ga_3\Ga_1\1)=
 \Ga_2(\Ga_1\Ga_3\Ga_1\1)\Ga_2
\eqtag\label{eq.4a2}
$$
given by the monodromy along
the boundary of the triangle~$S$ shown
in Figure~\ref{fig.4a2}. The monodromy about the
remaining singular fiber contained in the remaining
triangular region of the skeleton can be ignored due to
Proposition~\ref{presentation}.
Thus, the group $\pi_1(\Cp2\sminus B)$ is given by
$$
 \bigl\<\Ga_1,\Ga_2,\Ga_3\bigm|
 \text{\eqref{eq.e7.1}, \eqref{eq.e7.2}, \eqref{eq.4a2},
 \eqref{eq.leaf} with $m=n=2$}\bigr\>.
\eqtag\label{eq.G}
$$
Denote this group by~$G$. Using \GAP~\cite{GAP},
see Diagram~\ref{GAP.no.1}, one can see
that:
\roster
\item\local{G.1}
one has
$\ord G=41040=2^4\cdot3^3\cdot5\cdot19$, and the commutant $[G,G]$
is a perfect group
of order $6840=2^3\cdot3^2\cdot5\cdot19$;
\item\local{G.2}
the only perfect group of order $6840=\ord[G,G]$ is
$\SL(2,\F_{19})$;
\item\local{G.3}
relation \eqref{eq.4a2} in~\eqref{eq.G} follows from the others (as
dropping this relation does not change the
order of the group);
\item\local{G.4}
the order of each generator $\Ga_i$, $i=1,2,3$, in~$G$ equals
$6\cdot19$;
\item\local{G.5}
the group~$G$ is generated by $\Ga_1$ and~$\Ga_2$ only, as well as
by $\Ga_1$ and~$\Ga_3$ only.
\endroster
Due to~\loccit{G.3}, one can drop relation~\eqref{eq.4a2}; then
\eqref{eq.G} turns into the presentation given in the statement of
Theorem~\ref{th.e7.group}. The abelianization of~$G$ is the
cyclic group $\CG6$ generated by the image~$\bGa_1$ of~$\Ga_1$.
Hence, due to~\loccit{G.4},
the map $\bGa_1\mapsto\Ga_1^{19}$ splits the exact sequence
$$
\{1\}@>>>[G,G]@>>>G@>>>\CG6@>>>\{1\},
$$
representing~$G$ as a semi-direct product of its abelianization
and its commutant.

This completes the proof of Theorem~\ref{th.e7.group}.
\qed

\midinsert
\beginGAP%
gap> f := FreeGroup(3);;
gap> r1 := f.2*f.3*f.2^-1*f.3^-1;;
gap> r2 := (f.1*f.2*f.3)^3/(f.2*f.3^2);;
gap> r3 := f.1*(f.1*f.2*f.3)^3*f.1^-1*(f.1*f.2*f.3)^-3;;
gap> r4 := (f.1*f.2)^2*f.1/(f.2*(f.1*f.2)^2);;
gap> r5 := (f.1*f.3)^2*f.1/(f.3*(f.1*f.3)^2);;
gap> r6 := f.1*f.3*f.1^-1*f.2*f.1*f.3*f.1^-1/(f.2*f.1*f.3*f.1^-1*f.2);;
gap> g := f / [r1, r2, r3, r4, r5, r6];;
gap> Size(g);
41040
gap> List(DerivedSeriesOfGroup(g), AbelianInvariants); ## [G,G] is perfect
[ [ 2, 3 ], [  ] ]
gap> List(DerivedSeriesOfGroup(g), Size);              ## order of [G,G]
[ 41040, 6840 ]
gap> NumberPerfectGroups(6840); PerfectGroup(6840);    ## [G,G]=SL(2,19)
1
L2(19) 2^1 = SL(2,19)
gap> Size(f / [r1, r2, r3, r4, r5]);     ## drop r6
41040
gap> Size(Subgroup(g, [g.1]));           ## order of g.1
114
gap> Index(g, Subgroup(g, [g.1, g.2]));  ## g.1, g.2 generate G
1
gap> Index(g, Subgroup(g, [g.1, g.3]));  ## g.1, g.3 generate G
1
\endGAP
\figure[Diagram]
The \GAP\ output for $\bE_7\splus2\bA_4\splus2\bA_2$
\endfigure\label{GAP.no.1}
\endinsert

\proposition\label{incl.onto}
Let~$B$ be an irreducible plane sextic
with the set of singularities
$\bE_7\splus2\bA_4\splus2\bA_2$, and let $\MB$ be a Milnor ball
about the type~$\bE_7$ singular point of~$B$. Then the inclusion
homomorphism $\pi_1(\MB\sminus B)\to\pi_1(\Cp2\sminus B)$ is onto.
\endproposition

\proof
The statement follows from~\itemref{s.no.1}{G.5} above and
from~\eqref{eq.incl.e7}, which implies that $\Ga_1$ and~$\Ga_3$ are
in the image of the inclusion homomorphism.
\endproof

%
%

\section{Perturbations}\label{S.perturbations}

We compute the fundamental groups of all perturbations of a
type~$\bE_7$ singular point and apply these results to prove
Theorem~\ref{th.e7.pert}.

\subsection{Perturbations of a type~$\bE_7$ singularity}\label{s.e7.pert}
Let~$P$ be a type~$\bE_7$ singular point of a plane curve~$B$,
let~$\MB$ be a Milnor ball about~$P$, and let $B_t$, $t\in[0,1]$,
be a small perturbation of $B=B_0$ which remains transversal to
the boundary $\partial\MB$. We are interested in the perturbation
epimorphism $\pi_1(\MB\sminus B)\onto\pi_1(\MB\sminus B_t)$,
$t\ne0$.

According to E.~Looijenga~\cite{Looijenga},
the deformation classes (in the obvious sense) of perturbations of
a simple singular point~$P$ can be enumerated by the induced subgraphs
of the Dynkin diagram of~$P$ (up to a certain equivalence, which
is not important here). Comparing two independent classifications,
one can see that
any perturbation~$B_t$ of a type~$\bE_7$ singular point can be
realized by a family $C_t\subset\C^2$ of affine quartics
inflection tangent to the line at infinity (see,
\eg,~\cite{quintics}),
so that $(\MB,B_t)\cong(\C^2,C_t)$ for
each $t\in[0,1]$.
The groups $\pi_1(\C^2\sminus C_t)$ for such
quartics are found in~\cite{groups}, and all but five of them are
abelian.
The sets of singularities (perturbations) with nonabelian
fundamental group are those listed in Figure~\ref{fig.pert} and
$\bA_2\splus3\bA_1$, which is a further perturbation of
$\bD_5\splus\bA_1$ not changing the group.

\midinsert
\centerline{\vbox{\halign{\hss#\hss&&\qquad\hss#\hss\cr
\cpic{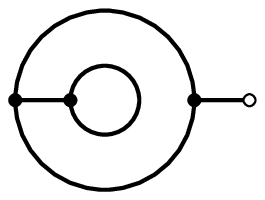}&\cpic{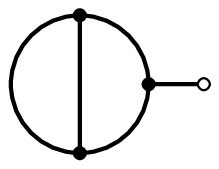}&
 \cpic{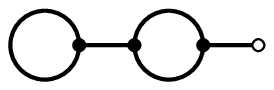}&\cpic{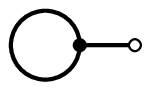}\cr
\noalign{\medskip}
$\bA_4\splus\bA_2$&$\bA_3\splus\bA_2\splus\bA_1$&
 $\bA_5\splus\bA_1$&$\bD_5\splus\bA_1$\cr
\crcr}}}
\figure\label{fig.pert}
Perturbations of $\bE_7$ with nonabelian fundamental group
\endfigure
\endinsert

In order to (re-)compute the four nonabelian groups, we realize
each perturbation by a family $\B_t\subset\Sigma_2$, $t\in[0,1]$,
of trigonal curves with a common type~$\tA{_1^*}$ singular
fiber~$F$, so that $\B_0$ is the isotrivial curve given by
${y'}^3=x'y'$. (We use the same coordinates $(x',y')$ as
in~\eqref{eq.coord}; in fact, $\B_t$ can be obtained from the
family~$C_t$ above by a birational transformation, similar
to~\ref{proof.1-1.e8}.) Then
$\MB\sminus B_t\cong\Sigma_2\sminus(\B_t\cup E\cup F)$,
and the group $\pi_1(\Sigma_2\sminus(\B_t\cup E\cup F))$ can be
found using the techniques outlined in~\ref{s.strategy}, without
adding the relation at infinity.

The skeleton of a typical curve~$\B_t$, $t\ne0$, which is maximal,
is one of those shown in Figure~\ref{fig.pert}. Let~$u$ be the
\white-vertex representing the fiber at infinity~$F$, and let~$v$
be the \black-vertex connected to~$u$. We take~$F_v$ for the
reference
fiber, and choose a
canonical
basis $\{\bb_1,\bb_2,\bb_3\}$
for~$\pi_v$ defined by the marking in which $[v,u]$ is the
edge~$e_ 1$ at~$v$. The braid monodromy is computed
using~\ref{s.monodromy}. In each case, there are two
relations~\eqref{eq.bc} given by the \cross-vertices connected
to~$v$ by the edges~$e'_1$ and~$e'_3$. For the third
\cross-vertex, if present, one should connect it to~$v$ by a solid
edge and a chain of bold edges and use~\eqref{eq.bb}.

The relations obtained are:
\Dashes
\dash
$\bA_4\splus\bA_2$:\enspace
 $\bb_1\bb_2\bb_1=\bb_2\bb_1\bb_2$,
 $(\bb_2\bb_3)^2\bb_2=\bb_3(\bb_2\bb_3)^2$,
 $\bb_2=\bb_3\bb_1\bb_3\1$;
\dash
$\bA_3\splus\bA_2\splus\bA_1$:\enspace
 $[\bb_1,\bb_3]=1$, $(\bb_1\bb_2)^2=(\bb_2\bb_1)^2$,
 $\bb_2\bb_3\bb_2=\bb_3\bb_2\bb_3$;
\dash
$\bA_5\splus\bA_1$:\enspace
 $[\bb_2,\bb_3]=1$, $(\bb_1\bb_2)^3=(\bb_2\bb_1)^3$,
 $\bb_3=\bb_1\bb_2\bb_1\1$;
\dash
$\bD_5\splus\bA_1$ and $\bA_2\splus3\bA_1$:\enspace
 $[\bb_1,\bb_2]=[\bb_1,\bb_3]=1$,
 $\bb_2\bb_3\bb_2=\bb_3\bb_2\bb_3$.
\endDashes

\Remark
The first group is shown in~\cite{groups} to be isomorphic to
$\Z\times\SL(2,\F_5)$. The last group is obviously
$\Z\times\BG3$. For the group~$G$ corresponding to the
set of singularities $\bA_5\splus\bA_1$,
one can easily deduce that
$\bb_1^3$ is a central element and then obtain a short exact
sequence
$$
\{1\}\longrightarrow\Z[t]/(t^3-1)\longrightarrow G
 \longrightarrow\Z\longrightarrow\{1\},
$$
the generator of the quotient~$\Z$ acting on the kernel \via\ the
multiplication by~$t$. (As a $\Z[t]$-module, the kernel is
generated by~$\bb_2$.) This result also agrees with~\cite{groups}.
\endRemark

\subsection{Proof of Theorem~\ref{th.e7.pert}}\label{proof.e7.pert}
It suffices to study the perturbations of the only
sextic~$B$ with nonabelian fundamental group, \ie, the one
considered in~\ref{s.no.1}. The set of singularities of~$B$ is
$\bE_7\splus2\bA_4\splus2\bA_2$.

The perturbations of the $\bA$-type points can be treated within
the framework of the paper, using the trigonal model~$\B$ of~$B$.
If one of the two type~$\bA_4$ points is perturbed, a
monovalent \cross-vertex appears in one of the two regions
adjacent to the insertion, see Figure~\ref{fig.4a2};
hence, the group
gets an additional relation $\Ga_1=\Ga_2$ or $\Ga_1=\Ga_3$ (\ie,
\eqref{eq.leaf} with $m=0$ or $n=0$) and, in view
of the commutativity relation $[\Ga_2,\Ga_3]=1$, becomes abelian.
The two cusps of~$\B$ are interchanged by the
complex conjugation, and it suffices to perturb
the one represented by the region~$S$ in Figure~\ref{fig.4a2}.
Then, a monovalent \cross-vertex appears in~$S$ and the
fundamental group
gets an extra relation $\Ga_1\Ga_3\Ga_1\1=\Ga_2$ (instead
of~\eqref{eq.4a2}). Thus, the group is a quotient of the group~$G$
given by~\eqref{eq.+loop} and hence is abelian, see~\ref{s.+loop}.

Now, consider a perturbation~$B'$ of the type~$\bE_7$ point~$P$.
Let~$\MB$ be a Milnor ball about~$P$. If the group
$\pi_1(\MB\sminus B')$ is abelian, then
$\pi_1(\Cp2\sminus B')$ is also abelian due to
Proposition~\ref{incl.onto}. For the four nonabelian groups
$\pi_1(\MB\sminus B')$, see~\ref{s.e7.pert}, we use the
description of the inclusion homomorphism found in~\ref{incl.e7}.
The `linear' basis $\{\bc_1,\bc_2,\bc_3\}$ considered
in~\ref{incl.e7} differs from the canonical basis
$\{\bb_1,\bb_2,\bb_3\}$ used in~\ref{s.e7.pert} by `half' the
monodromy about a type~$\bE_7$ singular point: one has
$$
\bb_1=(\bc_1\bc_2\bc_3)\bc_2(\bc_1\bc_2\bc_3)\1,\quad
\bb_2=(\bc_1\bc_2)\bc_3(\bc_1\bc_2)\1,\quad
\bb_3=\bc_1.
$$
Hence, in terms of the $\bb_i$, the inclusion homomorphism is
given by
$$
\bb_1\mapsto\Ga_1\Ga_2\1\Ga_1\Ga_3\Ga_1\1\Ga_2\Ga_1\1,\quad
\bb_2\mapsto\Ga_1\Ga_2\1\Ga_1\Ga_2\Ga_1\1,\quad
\bb_3\mapsto\Ga_1.
$$
(We used the commutativity relations~\eqref{eq.e7.1} to simplify
the expressions obtained.)
Now, it remains to express the extra relations in terms of the
$\Ga$-basis, add them to the presentation given by~\eqref{eq.G},
and use \GAP~\cite{GAP}. (In fact, in each of the four cases
listed in~\ref{s.e7.pert}, the
first extra relation alone make the group abelian.)
\qed

\widestnumber\key{EO1}
\refstyle{C}

\enddocument